\documentclass{amsart}

\usepackage{fullpage,url,xcolor,hyperref}
\usepackage{algpseudocode}
\usepackage{algorithm,mathtools}
\usepackage{multirow}
\usepackage{subcaption}
\usepackage[section]{placeins} 

\algdef{SE}[DOWHILE]{Do}{doWhile}{\algorithmicdo}[1]{\algorithmicwhile\ #1}%

\DeclareMathOperator{\dgm}{dgm}
\newcommand{\Hg}{\mathrm{H}}
 \usepackage{amsmath,amsthm,bbm,amssymb}
\DeclareMathOperator{\birth}{birth}
\DeclareMathOperator{\death}{death}

\DeclareMathOperator{\logg}{\log\log}

\def\R{\mathbb{R}}

\def\P{\mathbb{P}}

\def\S{\mathbb{S}}

\def\cD{\mathcal{D}}
\def\cF{\mathcal{F}}

\def\cI{\mathcal{I}}

\def\cL{\mathcal{L}}

\def\cT{\mathcal{T}}

\def\cU{\mathcal{U}}
\def\cS{\mathcal{S}}

\def\cX{\mathcal{X}}

\def\cL{\mathcal{L}}

\def\mS{\mathrm{S}}
\def\mN{\mathrm{N}}

\def\b{\mathrm{b}}
\def\d{\mathrm{d}}

\newcommand{\given}{\;|\;}

\newcommand{\bX}{{\mathbf{X}}}

\newcommand{\set}[1]{\left\{#1\right\}}

\newcommand{\norm}[1]{\left\|#1\right\|}
\newcommand{\param}[1]{\left(#1\right)}

\newcommand{\cprob}[2]{\mathbb{P}\left(#1\given #2\right)} 

\newcommand{\bz}{\mathbf{z}}
\newcommand{\by}{\mathbf{y}}
\newcommand{\bx}{\mathbf{x}}
\newcommand{\bv}{\mathbf{v}}

\providecommand{\setthms}[1]{#1}
\setthms{
\newtheorem{lem}{Lemma}[section]

\newtheorem{con}[lem]{Conjecture}

\theoremstyle{definition}

}
\newcommand{\cech}{\v{C}ech }

\newcommand{\ninf}{n\to\infty}

\newcommand{\limninf}{\lim_{\ninf}}
\definecolor{mygreen}{rgb}{0, 0.68, 0.31}
\definecolor{myred}{rgb}{1.0, 0,0}

\numberwithin{equation}{section}

\def\bsplit#1\esplit{\begin{split} #1 \end{split} }
\def\splitb#1\splite{\begin{split} #1 \end{split} }
\def\beq#1\eeq{\begin{equation} #1 \end{equation}}
\def\eqb#1\eqe{\begin{equation} #1 \end{equation}}

\newtheorem*{remark}{Remark}

\usepackage[foot]{amsaddr}

\title{On the universality of random persistence diagrams}

\author{Omer Bobrowski$^{1,2}$}
\address{$^1$ School of Mathematical Sciences, Queen Mary University of London.}
\address{$^2$ Viterbi Faculty of Electrical and Computer Engineering, Technion -- Israel Institute of Techonology.}
\author{Primoz Skraba$^1$}
\email{o.bobrowski@qmul.ac.uk, p.skraba@qmul.ac.uk}

\begin{document}
\maketitle
\begin{abstract}
    One of the most elusive challenges within the area of topological data analysis is understanding the distribution of persistence diagrams. Despite much effort, this is still largely an open problem. In this paper, we present a series of novel conjectures regarding the behaviour of persistence diagrams arising from random point-clouds. We claim that these diagrams obey a universal probability law, and include an explicit expression as a candidate for what this law is. We back these conjectures with an exhaustive set of experiments, including both simulated and real data. We demonstrate the power of these conjectures by proposing a new hypothesis testing framework for individual features within persistence diagrams.
\end{abstract}

\vspace{10pt}

Topological Data Analysis (TDA) focuses on extracting structural topological information from data, in order to enhance their processing in statistics and machine learning.
This field has been rapidly developing over the past two decades, bringing together mathematicians, statisticians, computer scientists, engineers, and data scientists. 
The motivation for incorporating topological features in data analysis  is that topological methods are highly versatile, coordinate-free, and robust to various deformations. Topological methods have been applied successfully in numerous applications developed over the last decade, in areas such as neuroscience ~\cite{bardin2019topological,gardner2022toroidal,giusti2015clique,petri_homological_2014}, medicine and biology ~\cite{chan2013topology,ichinomiya2020protein,mcguirl2020topological,vipond2021multiparameter},  material sciences ~\cite{herring2019topological,kramar2013persistence,lee2017quantifying,onodera2019understanding}, dynamical systems \cite{kramar2016analysis,tymochko2020using,xian2020capturing}, and cosmology ~\cite{adler_modeling_2017, pranav_topology_2016, van_de_weygaert_alpha_2011}.

In computing topological features for a given dataset, one of the key challenges in TDA is how to distinguish between signal and noise. By ``signal,'' we broadly refer to features that represent meaningful structures underlying the data, while ``noise'' refers to features that arise from the local randomness and inaccuracies within the data. The identification of this structure has been a central problem in TDA (for example, see \cite{chazal_sampling_2009,chazal_persistence-based_2013,niyogi_finding_2008, stolz2020geometric}).
Over the years, statisticians and probabilists have studied various approaches to address this challenge.
The main difficulty is that
the distribution of topological noise seems to have many different shapes and forms, and is highly sensitive to how data are
generated.  Our main goal in this paper is to \emph{refute} the last premise, arguing that the distribution of noise in persistent homology of geometric complexes, viewed in the right way, is in fact \emph{universal}. 

Given a point-cloud $\cX$, \emph{persistent homology} searches for structures such as holes and cavities formed by $\cX$, and records the scale at which they are created and terminated (referred to as the \emph{birth} and \emph{death} scales, respectively).
The output is known as a \emph{persistence diagram} -- a collection of points $p=(\b,\d)\in \R^2$ marking the {\bf b}irth and {\bf d}eath times of these structures. The prevalent intuition in TDA is that the further a point $p$ is from the diagonal $\d=\b$, the more meaningful the corresponding topological feature is. 
Turning this intuition into a rigorous and applicable statistical methodology has been the focus of many studies in TDA (cf.~\cite{blumberg_robust_2013,fasy_confidence_2014,reani2021cycle, vejdemo2020multiple}). However, two decades of research have yet to provide a generic, robust, and theoretically justified framework. Another, highly related, line of research has been the theoretical probabilistic analysis of persistence diagrams generated by random data, as means to establish a null-distribution for persistent homology. While developing the mathematical theory has been  fruitful (cf.~\cite{adler_crackle:_2014,bobrowski_maximally_2017, hiraoka_limit_2018,owada_limit_2017,yogeshwaran_topology_2015,yogeshwaran_random_2016}), its use in practice has been  limited. The main gap between  theory and practice is that all of these studies indicate that the distribution of noise in persistence diagrams: (a) does not have a simple closed form description, and (b) strongly depends on the model generating the point-cloud. 

Our main contribution in this paper is as follows. Given a persistence diagram generated by 
a random point-cloud, consider only the noisy features, and measure their \emph{persistence} using the  death/birth ratio. We argue that  taking the limit (as the point-cloud gets larger), the distribution of these persistence values is \emph{universal}, in the sense that it is independent of the model generating the point-cloud. This result is loosely analogous to the sums of many different types of random variables always converging to the normal distribution, via the central limit theorem.

The emergence of universality in this context is highly surprising and unexpected. In this paper, we phrase universality as a sequence of conjectures, which are  supported by an extensive experimental study. The proofs for these conjectures will require significant advance and a shift of paradigm in the study of stochastic topology, laying the groundwork for many potential future discoveries.

Recall that one of the fundamental uses of the central limit theorems in statistics, is for hypothesis testing. Based on our conjectures, we develop a similar hypothesis testing framework for persistence diagrams, to test whether individual features should be classified as signal or noise. We believe the framework we provide here is very powerful, as it allows to compute a numerical significance measure (p-value) for each individual feature in a persistence diagram, using very few assumptions on the underlying model.

\section{Persistent homology for point-clouds}

A \emph{point-cloud} is a finite collection of points $\cX$ in a metric space. For simplicity, in this paper we will always use point-clouds in the Euclidean space $\R^D$, which is the most common setting in applications.
In order to turn point-clouds into shapes that can be studied topologically, the common practice in TDA is to  use their geometry to generate simplicial complexes.

An \emph{abstract simplicial complex} 
can be thought of as a ``high-dimensional graph'' where in addition to vertices and edges, we include triangles, tetrahedra, and higher dimensional simplices. The only restriction is that it is closed under inclusion, i.e., for every simplex, the complex must include all its lower dimensional faces as well.

Given a point-cloud $\cX$, we can construct a simplicial complex whose vertex set is $\cX$, and its simplices are determined by the geometric configuration of $\cX$.
The two most common constructions in TDA are both parameterized by a radius value $r>0$, and defined as follows:
\begin{itemize}
    \item The \textbf{\v Cech complex}: 
    A subset of $(k+1)$ points spans a $k$-simplex  if the intersection of all $r$-balls around them is non-empty.
    \item The \textbf{Vietoris-Rips complex}: 
    A subset of $(k+1)$ points spans a $k$-simplex  if all \emph{pairwise} intersections of the $r$-balls around them are  non-empty.
\end{itemize}

Given a simplicial complex $X$, one can study its \emph{homology groups}, denoted $\Hg_k(X)$. Loosely speaking, homology is an algebraic object that represents structural information about $X$. For example, $\Hg_0(X)$ contains information about connected components, $\Hg_1(X)$ -- about closed loops surrounding holes, $\Hg_2(X)$ -- about closed surfaces surrounding cavities. Generally, we say that $\Hg_k(X)$ represents information about `$k$-cycles'. For additional background on algebraic topology, see \cite{hatcher_algebraic_2002, munkres_elements_1984}. In this paper, we only consider $k\ge 1$, i.e., everything \emph{but} connected components. 

\begin{figure}
    \centering
    \includegraphics[width=0.7\linewidth]{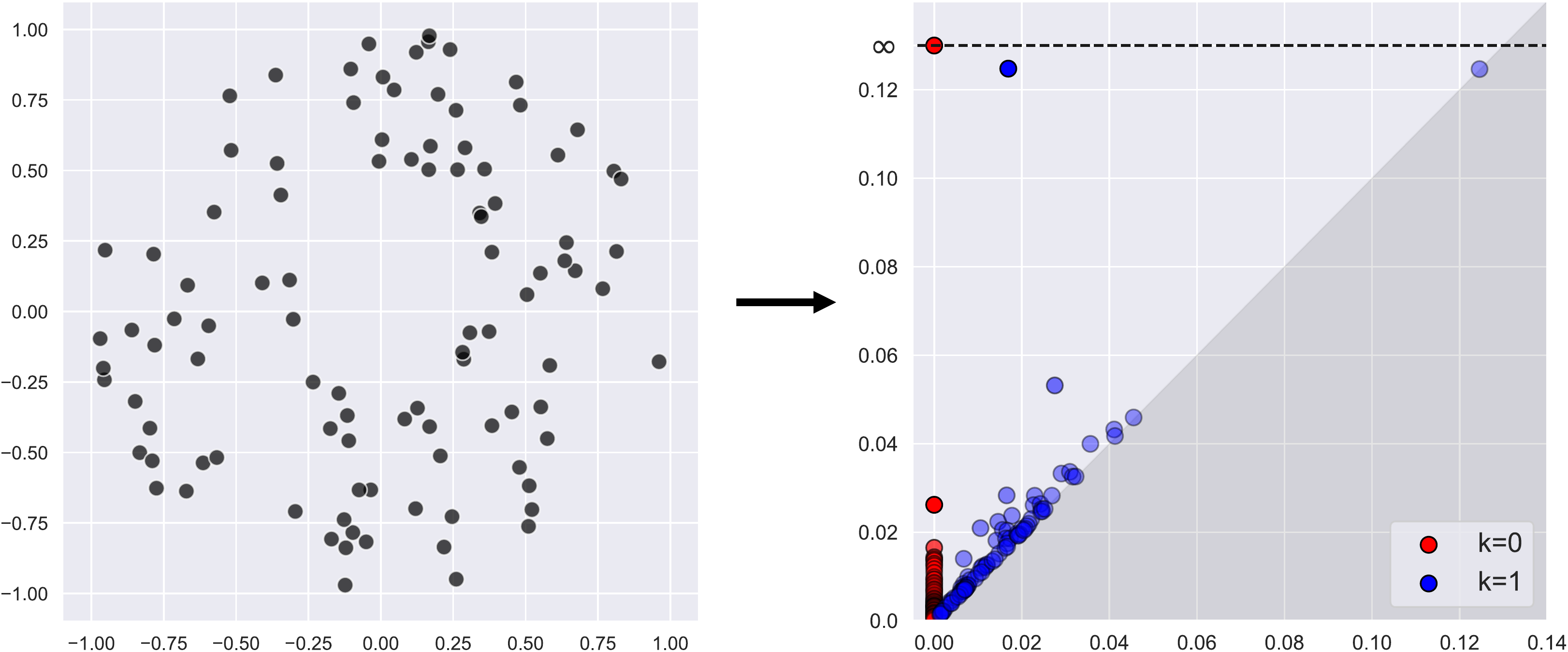}
    \caption{A point-cloud of 100 points, sampled from an annulus, and the persistence diagram corresponding to the \v Cech filtration. The red points mark persistent $0$-cycles (components) and the blue points mark persistent $1$-cycles (holes). A single red and single blue point stand out, reflecting the structure of the annulus.}
    \label{fig:PD}
\end{figure}

In the context of the \v Cech and Rips complexes, the homology $\Hg_k(X)$ is highly sensitive to the choice of the radius parameter $r$. To overcome this sensitivity, the approach proposed by \emph{persistent homology} is to consider the entire range of parameter values, and track the evolution of $k$-cycles as the value of $r$ increases. In this process (called a `filtration'), cycles are created (born) and later are filled-in (die). Once persistent homology is calculated, it is most commonly presented via \emph{persistence diagrams}. These are 2-dimensional scatter plots, where each point is of the form $p=(\b,\d)\in \R^2$, representing the \emph{birth} and \emph{death} times (radii, in our case) of a single cycle in the filtration. 
See Figure \ref{fig:PD} for an example.
For a given filtration $\cF$, we denote its $k$-th persistence diagram by $\dgm_k(\cF)$. For more details see \cite{edelsbrunner_computational_2010,zomorodian_computing_2005}.

One of the original motivations for using persistent homology is that it allows us to detect meaningful structures that emerge in the data. The simplest way to do so, is by looking for points $p\in\dgm_k$ that are far away from the diagonal. These points represent cycles with a long lifetime ($\d-\b$), which are intuitively ``significant'', compared to points near the diagonal that are created due to the finite and noisy nature of a dataset (see Figure \ref{fig:PD}). While this approach is very intuitive,
justifying it theoretically and providing \emph{quantitative statements}, are among the greatest challenges in the field.

\section{Noise distribution in persistence diagrams}
The problem of detecting the meaningful features in persistence diagrams can be formalized as follows.
For any given persistence diagram $\dgm_k$, we assume a decomposition into $\dgm_k^{_\mS} \cup\dgm_k^{_\mN}$. The points in $\dgm_k^{_\mS}$ correspond to the ``signal'' -- the meaningful structure that we wish to extract from the data. The points in $\dgm_k^{_\mN}$ are the ``noise'' -- cycles that are created due to the randomness in the data, and contain no useful information. The challenge is  to decide for each $p\in\dgm_k$ whether $p\in\dgm_k^{_\mS}$ or $p\in\dgm_k^{_\mN}$, and to provide quantitative statistical guarantees on this decision.
Revealing the distribution of $\dgm_k^{_\mN}$ would serve as a null-distribution for persistence diagrams, that will allow us to accurately detect the signal. 

\subsection{Previous Work}

Closest in spirit to our current study is \cite{hiraoka_limit_2018}, where  persistence diagrams were considered as random measures. The main result shows that for a wide class of stationary processes there exists a non-random limiting measure.
 Later \cite{divol2019density} proved that these measures have density, and \cite{shirai2021limit} extended the limit to marked point processes. In \cite{owada2020convergence}  the limit of persistence diagrams as compact sets was studied, in an extreme value setting. Other related works consider various marginals of persistence diagrams (e.g., Betti numbers, Euler characteristic) and study their limiting behavior in various settings \cite{bobrowski_distance_2014,bobrowski_maximally_2017,bobrowski_homological_2020,bobrowski2022homological,kahle_limit_2013,owada_limit_2017,yogeshwaran_topology_2015, yogeshwaran_random_2016}. While these results have established a  rich mathematical theory for persistence diagrams, translating them into statistical tools has lagged behind. The two main reasons are: (a) these results show that various limits exist, but in most cases without any explicit description, (b) these limits are highly dependent on the distribution generating the point-cloud.

Quite a different approach was taken in \cite{adler_modeling_2017}, where persistence diagrams were modelled as Gibbs measures, whose parameters can be estimated from the data. This framework is quite generic and powerful. Here as well, the resulting distribution varies between different point-clouds, and the parameters need to be estimated for every model separately.

With respect to topological inference, the TDA literature provides various powerful methods, all based on the premise that the distribution of persistence diagram is inaccessible.
In \cite{fasy_confidence_2014} the authors  introduced confidence intervals for persistence diagrams, overcoming this issue by using sub-sampling methods and the stability of the persistence diagrams \cite{cohen-steiner_stability_2007}. The confidence intervals are computed for an entire diagram, rather than for individual features, and thus are quite coarse. In addition, this approach requires restrictive assumptions on the underlying space (closed manifold), as well as  the probability distributions and its derivatives. 
Other useful methods include 
 distance to measures \cite{chazal_geometric_2011,chazal_robust_2017}, witness complexes \cite{de_silva_topological_2004}, persistence landscapes \cite{bubenik_statistical_2015},  multi-cover bifiltration \cite{corbet2021computing,sheehy2012multicover}, and other sub-sampling based methods \cite{blumberg_robust_2013,agerberg2022data,chazal_convergence_2015,chazal_robust_2017-1,solomon2021geometry,abrams_finding_2002}.

\section{Universal distribution of persistent cycles}

\begin{figure}
    \centering
    \includegraphics[width=0.7\columnwidth]{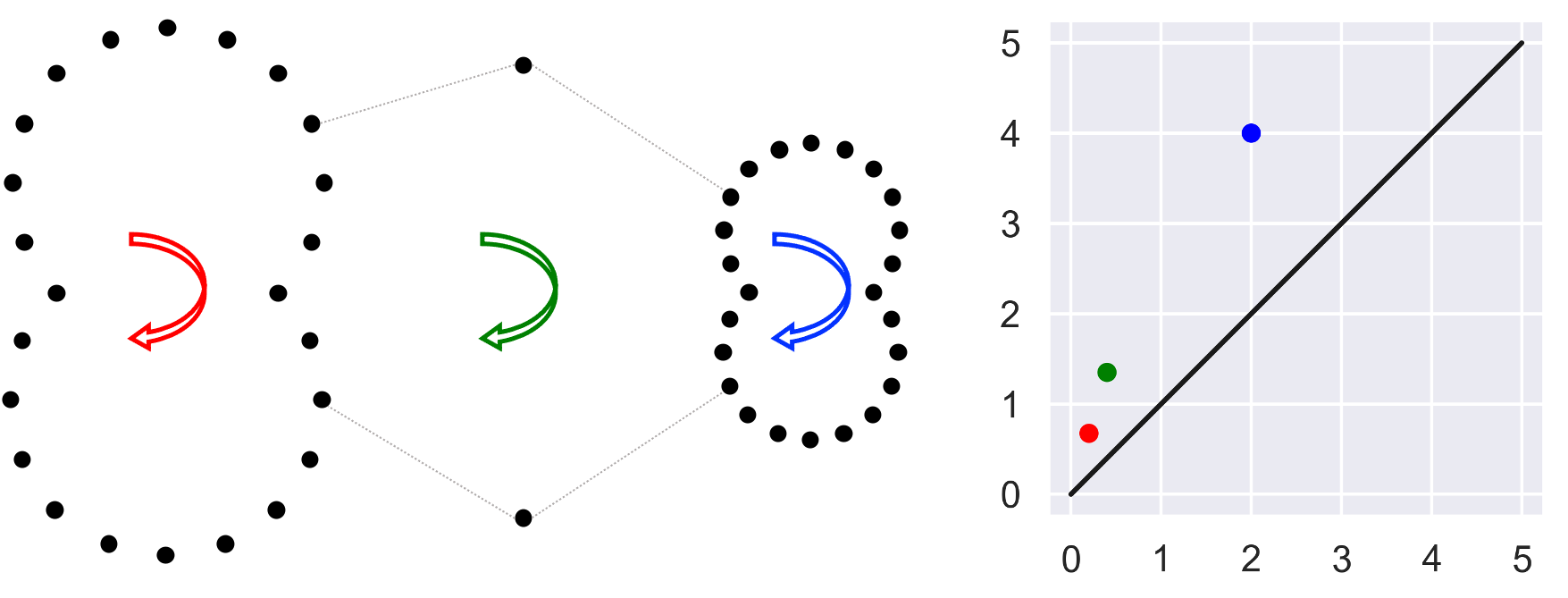}
    \caption{Multiplicative persistence. The point-cloud on the left generates three 1-dimensional cycles (marked by the arrows). On the right we see the corresponding persistence diagram with matching cycle colors. Observe that if we measure the lifetime of cycles, the green cycle will be considered the most significant, while it seems to be generated by outliers. In addition, the structure of the other two cycle (red and blue) is identical, yet their lifetime is different. Both issues are resolved if we take the death/birth ratio instead.}
    \label{fig:multip}
\end{figure}

As described above, a common practice in TDA is to measure the significance of topological features by considering the lifetimes of persistent cycles, i.e., $\Delta = (\death-\birth)$. This method is intuitive in toy examples (see Figure \ref{fig:PD}), as it captures the geometric ``size'' of topological features. However, a strong case can be made  \cite{bobrowski_maximally_2017} that for geometric complexes the ratio $\pi = (\death/\birth)$ is in fact a more robust statistic to discriminate between signal and noise in persistence diagrams (for $k>0$). There are two main justification for this statement. Firstly, the ratio $\pi$ is scale invariant, so that cycles that have exactly the same structure but exist at different scales are weighed the same. Secondly, datasets often contain outliers that may generate cycles with a large diameter, and consequently their lifetime $\Delta$ will also be large. However, the value of $\pi$ for such outliers should remain low, compared to features that occur in dense regions (see Figure \ref{fig:multip}).  

Let $\cS$ be a $d$-dimensional metric-measure space, and let $\bX_n = (X_1,\ldots,X_n)\in S^n$ 
a sequence of random variables (points), whose joint probability law is denoted by $\P_n$. Let $\P = (\P_n)_{n=1}^\infty$, and denote $\S=(\cS,\P)$, which we refer to as the \emph{sampling model}.
Fix a filtration type $\cT$ (e.g. \cech or Rips), and a homological degree $k>0$, and compute the $k$-th \emph{noise} persistence diagram $\dgm_k^{_\mathbf{N}}(\bX_n;\cT)$, which in short we denote by $\dgm_k$.
We study the distribution of the random persistence values $\set{\pi(p)}_{p\in \dgm_k}$ (where $\pi(p) = \death(p)/\birth(p)$) and refer to them as the $\pi$-values of the diagram.

In  \cite{bobrowski_maximally_2017}, the asymptotic scaling for the largest $\pi$-value was studied. Denoting $\pi_{k}^{\max} := \max_{p\in\dgm_k} \pi(p)$, the main result  (Theorem 3.1) shows that with high probability
\[
    A\param{\frac{\log n}{\logg n}}^{1/k} \le \pi_{k}^{\max} \le B \param{\frac{\log n}{\logg n}}^{1/k},
\]
for some constants $A,B>0$.
This is in contrast to \emph{signal} cycles, where using the results in \cite{bobrowski2022homological}, we know that they are of the order of $n^{1/d}$.
Thus, the $\pi$-values provide a strong separation (asymptotically) between signal and noise in persistence diagrams. 

\begin{remark} The results in \cite{bobrowski_maximally_2017} and \cite{bobrowski2022homological} are phrased for the uniform distribution on a $d$-dimensional box and a $d$-dimensional flat torus, respectively.
Nevertheless, the proofs can be adapted to a wide class of compact spaces and distributions, with the same asymptotic rates. 
\end{remark}

\subsection{Weak universality}\label{sec:weak}
We start by considering the case where $\P_n$ is a product measure, and the points $X_1,\ldots,X_n$ are iid (independent and identically distributed). Given $\dgm_k$ as defined above, denote the empirical measure of $\pi$-values as
\[
\Pi_{n} = \Pi_{n}(\S,\cT,k) := \frac{1}{|\dgm_k|}\sum_{p\in\dgm_k}\delta_{\pi(p)}.
\]

\begin{figure*}[h]
\centering
 		\includegraphics[width=\textwidth]{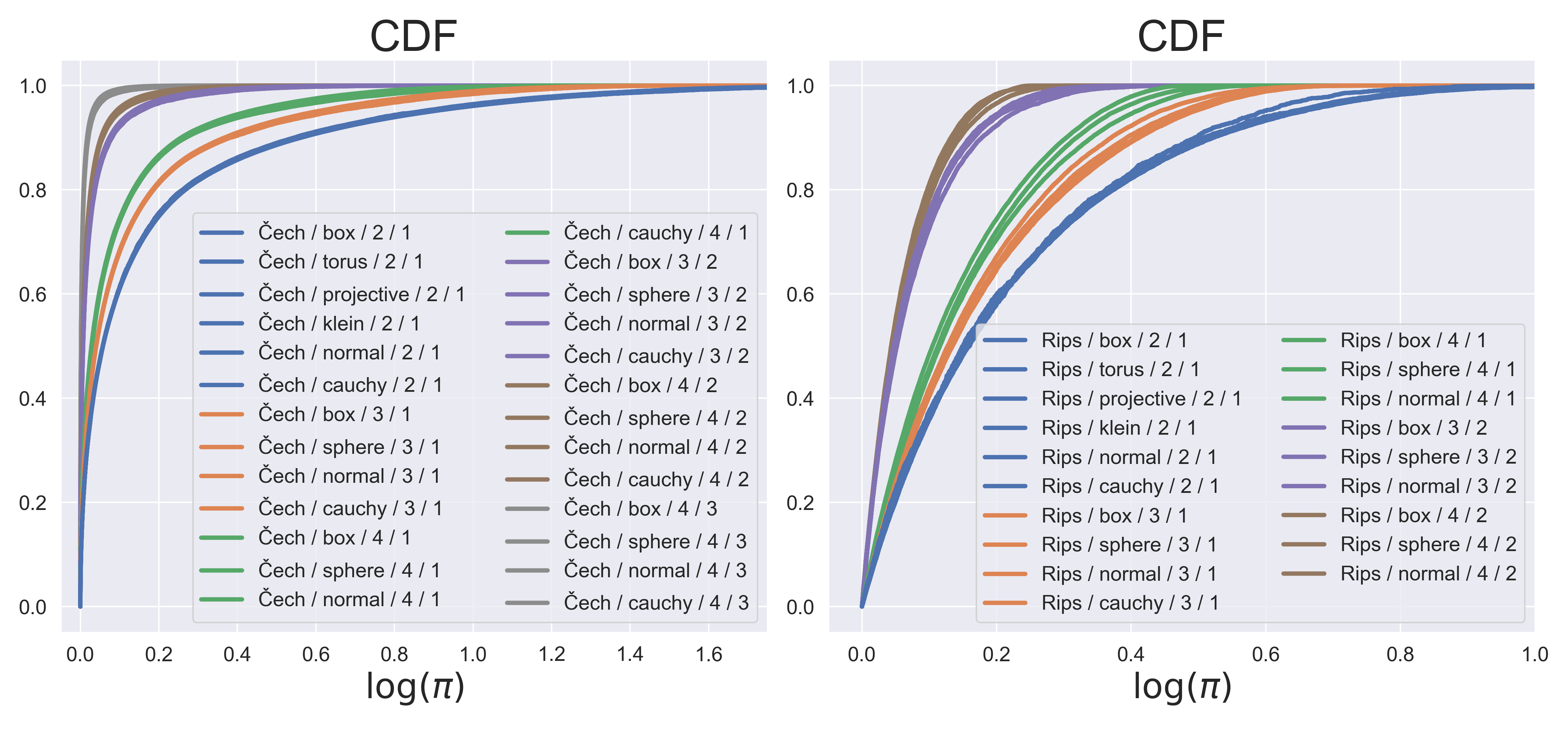}
 	\caption{The distribution of $\pi$-values. We take the empirical CDFs of the $\pi$-values (log-scale), computed from various  iid samples. The legend format is $\cT / \P / d / k$, where $\cT$ is the complex type, $\P$ is the probability distribution, $d$ is the dimension of the sampling space, and $k$ is the degree of homology computed.
 	  By `box',`torus', `sphere', `projective(-plane)', and `klein(-bottle)' we refer to the uniform distribution on the respective space or its most natural parametrization, while `normal' and `cauchy' refer to the corresponding non-uniform distributions.
 	  For more details about the experiments, and other models tested, see Appendix \ref{app:iid}. \label{fig:pi_cdfs}}
\end{figure*}

Note that $\Pi_n$ is a (discrete) probability measure on $\R$. In Figure \ref{fig:pi_cdfs} we present the CDF of $\Pi_{n}$, for various choices of $\S,\cT$, and $k$. In the figure we see that if we fix $d$ (dimension of $\cS$), $\cT$, and $k$, then the resulting CDF depends on neither the space $\cS$ nor the distribution $\P_n$. 
This observation leads to our first conjecture. 
Note that while our experiments indicate that this phenomenon spans various sampling models $\S$, we  obviously do not expect it to hold for \emph{every possible} sample. We thus denote by $\cI_d$ the class of relevant $d$-dimensional sampling models $\S$ (the extent of which is to be determined as future work).

\begin{con}\label{con:iid}
Fix $d,\cT,$ and $k>0$.  For any $\S\in \cI_d$, we have
\[
    \limninf\Pi_{n} = \Pi^*_{d,\cT,k}.
\]
\end{con}
The exact notion of limit is to be determined.
We further conjecture that the class $\cI_d$ should be fairly large. From our experiments, it seems that the space $\cS$ may belong to  a large class of spaces, including smooth Riemannian manifolds and stratified spaces (compact or otherwise). The joint distribution $\P_n$ should be continuous and iid, but otherwise is fairly generic (possibly even without any moments, see the Cauchy example in Figure \ref{fig:pi_cdfs}).

To conclude, we conjecture that a wide range of iid point-clouds, exhibits a {\bf universal} limit for  the $\pi$-values. We name this phenomenon ``weak universality'', since while the limit is independent of $\S$ (hence, universal), it does depend on $d,\cT,k$ as well as the iid assumption.
This is in contrast to the results in the next section.

\subsection{Strong universality}\label{sec:strong}
\begin{figure*}[h]
\begin{center}
	\begin{subfigure}{\textwidth}
	\centering
 		\includegraphics[width=0.8\textwidth]{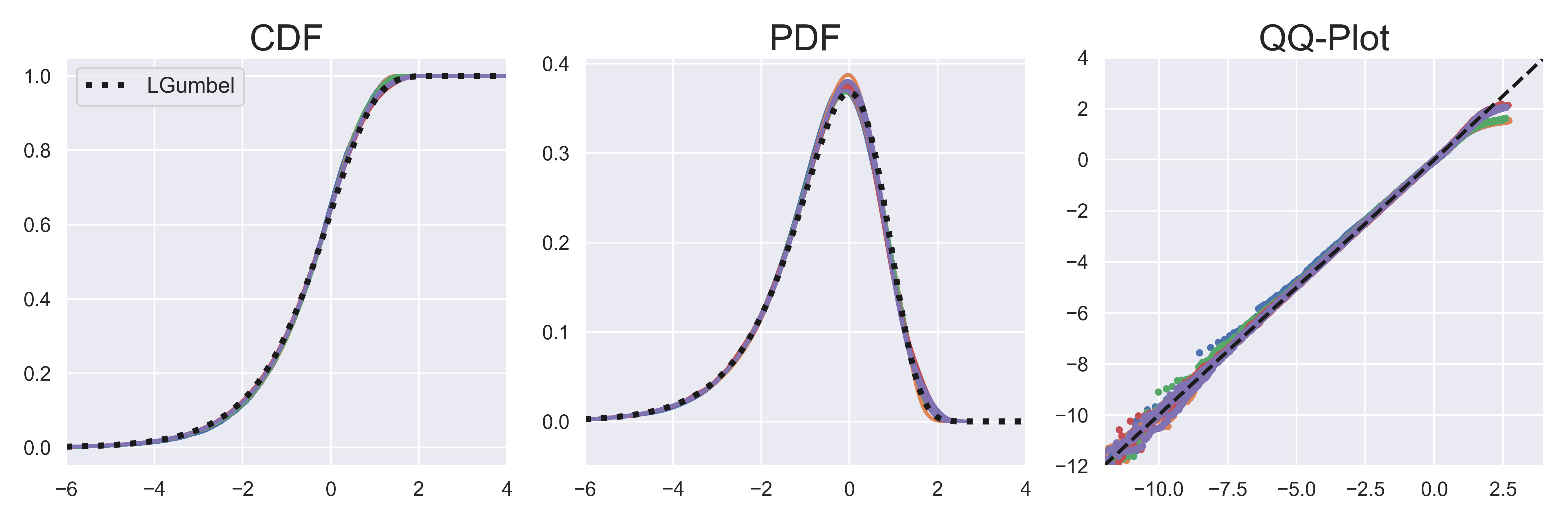}
\end{subfigure}
		\begin{subfigure}{\textwidth}
	\centering
		\includegraphics[width=0.8\textwidth]{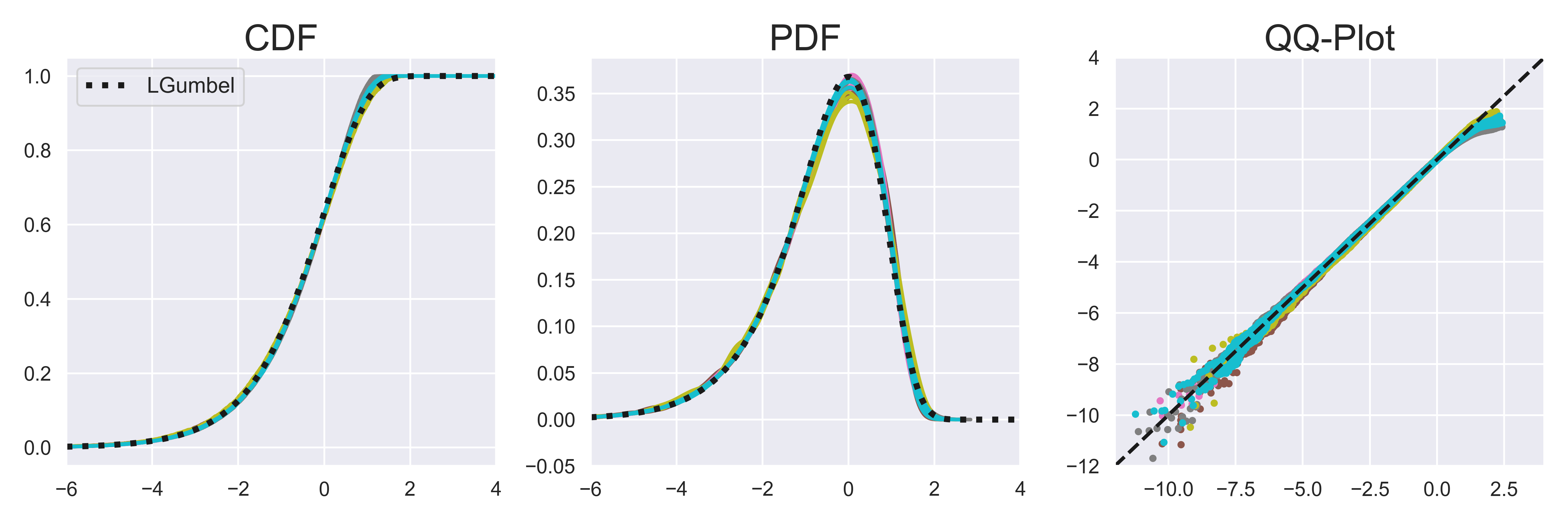}
	\end{subfigure}
	\end{center}
 	\caption{  The distribution of $\ell$-values. On the left are the empirical CDFs, and in the center are the kernel density estimates for the PDFs. Top row - \v Cech complex, bottom row - Rips complex. The curves contains \emph{all} the iid samples included in Figure \ref{fig:pi_cdfs} (i.e.~26 for the \v Cech and 20 for the Rips). The dashed line corresponds to the left-skewed Gumbel distribution. The figure on the right is the QQ plot for all the different samples when compared to the left-skewed Gumbel.\label{fig:loglog}}
\end{figure*}

The following procedure is not highly intuitive, and was discovered partly by chance. Nevertheless,  the results that this procedure yields are quite convincing and powerful.
Take a random persistence diagram $\dgm_k$ as described above. For each $p\in\dgm_k$ apply the following transformation
\eqb\label{eqn:loglog}
    \ell(p) := A\logg(\pi(p)) + B,
\eqe
where
\eqb\label{eqn:AB}
    A = \begin{cases}1 & \cT=\text{Rips},\\
    1/2 & \cT = \text{\v Cech},\end{cases}\qquad B = -\lambda - A\bar L
\eqe
and where $\lambda$ is the Euler-Mascheroni constant (=0.5772156649...), and $\bar L = \frac1{|\dgm_k|} \sum_{p\in \dgm_k}\logg(\pi(p))$.
We refer to the set $\set{\ell(p)}_{p\in \dgm_k}$ as the $\ell$-values of the diagram.
In Figure \ref{fig:loglog} we present the empirical CDFs   (left column) of the $\ell$-values, as well as the kernel density estimates for their PDFs (center column), for all the iid samples that were included in Figure \ref{fig:pi_cdfs}. We observe that all the different settings ($\S,\cT,k$) yield exactly the same distribution under the transformation given by \eqref{eqn:loglog}. We refer to this phenomenon as ``strong universality''.

While strong universality for iid point-clouds is by itself a very unexpected and useful behavior, we wish to investigate whether it applies in more general and realistic scenarios as well. We start by generating non-iid point-clouds that include: (a) sampling the path of a Brownian motion,  and (b) sampling the Lorenz dynamical system.
We then proceed to real data, and examine patches of natural images as well as time delay embeddings of audio recordings. For more details about these experiments, see Section \ref{sec:exp}.
While the distribution of the $\pi$-values for these models is vastly different than the one we observed in Section \ref{sec:weak}, 
all of these examples exhibit the same strong universality behavior. 
We present a subset of these results in Figure \ref{fig:noniid} and the full body of results can be found in Appendix \ref{app:exp}.

\begin{figure*}[h]
\begin{center}
	\begin{subfigure}{\textwidth}
	\centering
 		\includegraphics[width=\textwidth]{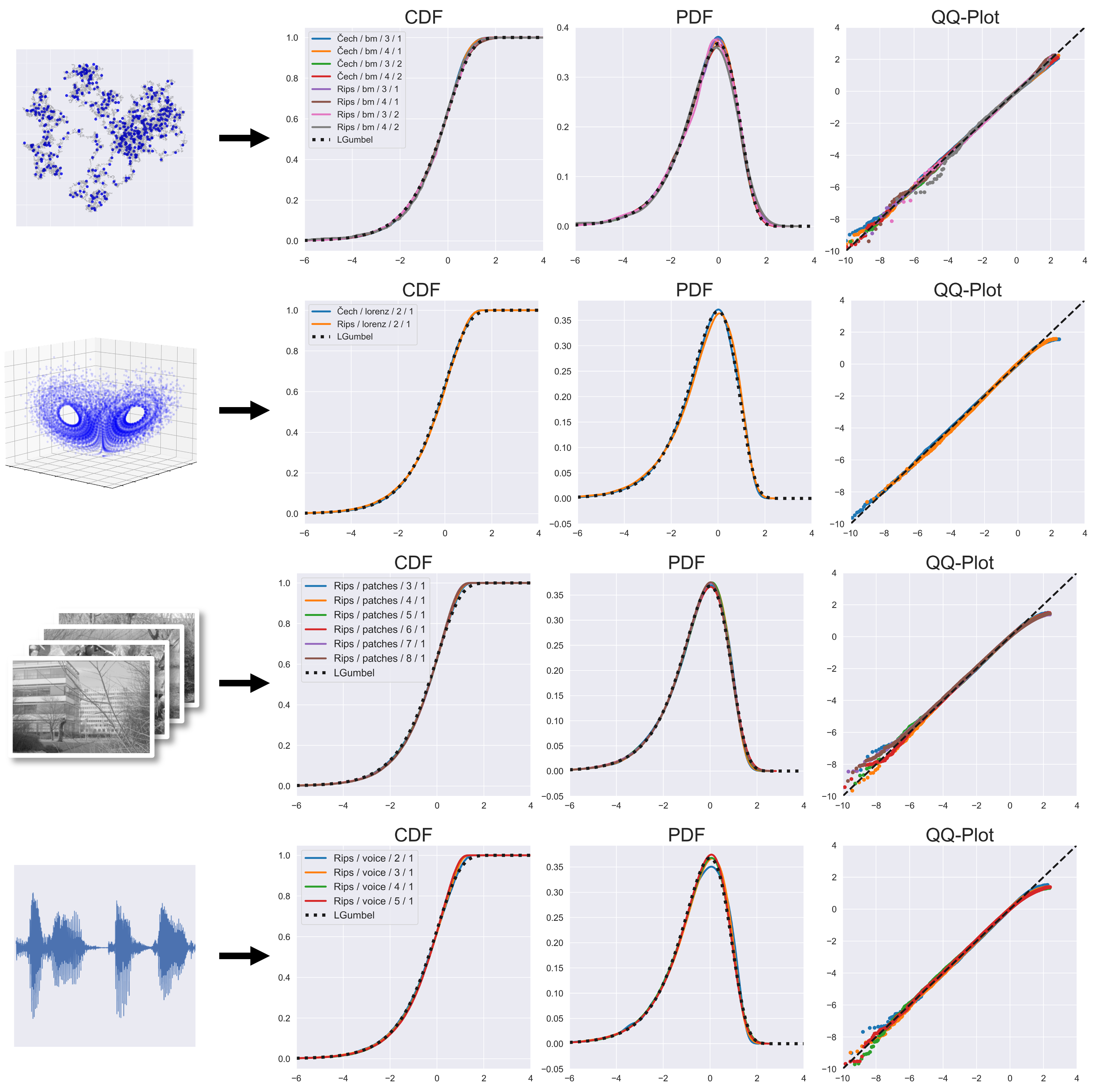}
\end{subfigure}
	\end{center}
 	\caption{Examples of non-iid point-clouds and real-data. First row - sampling the path of a $d$-dimensional Brownian motion. Second row - sampling from the Lorenz dynamical system.
 	Third row - sampling patches extracted from natural images. Fourth row - taking time-delay embeddings of audio recordings.
 	The figures on the left are visual examples of the data used. The three figures on the right are the corresponding estimates for the CDF, PDF and QQ-plot.
 	\label{fig:noniid}}
\end{figure*}

To conclude, our experiments highly indicate that persistent $\ell$-values have a universal limit for a wide class of sampling models $\S$. We shall denote this class by $\cU$. For our main conjecture, we consider the empirical measure of the $\ell$-values,
\[
\cL_n = \cL_n(\S,\cT,k) := \frac{1}{|\mathrm{dgm_k}|}\sum_{p\in\dgm_k}\delta_{\ell(p)}.
\]

\begin{con}\label{con:loglog}
For any $\S\in \cU$, and for any $k\ge1 $, and $\cT\in \{\text{\v Cech, Rips}\}$,
\[
\limninf \cL_n = \cL^*,
\]
where $\cL^*$ is independent of $\S$, $\cT$, and $k$.
\end{con}

\begin{figure}[h]
\begin{center}
		\includegraphics[width=0.5\columnwidth]{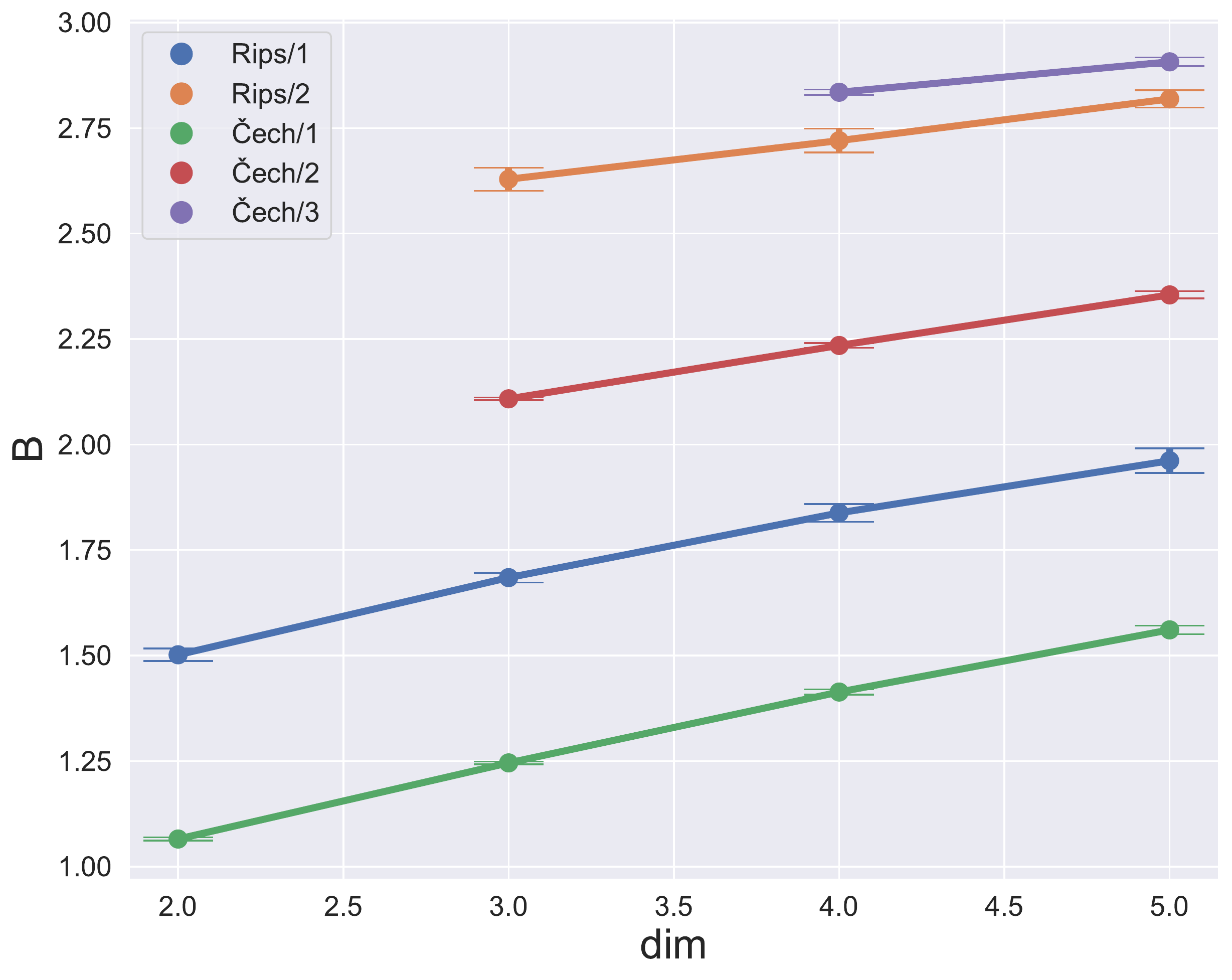}
 	\caption{The parameter value $B$ for iid samples with  various choices of $(d,\cT,k)$. We present the mean $\pm$standard deviation of our estimates across the various models. 
\label{fig:means}}
 	\end{center}
\end{figure}

Note that the only part in Conjecture \ref{con:loglog} that depends on the distribution generating the point-cloud is the value of  $B$ in \eqref{eqn:loglog} (compare this to the mean and variance in the classical CLT, for example).
It is therefore a natural question to examine
the dependency of $B$ in the model parameters $(\S, \cT,k)$. In  Figure~\ref{fig:means}, we examine the value of $B$ for different iid settings. As suggested by Conjecture \ref{con:iid}, the value of $B$ depends on $d,\cT,k$, but is otherwise independent of $\S$. Figure \ref{fig:means} also suggests a linear relationship between $B$ and the sampling dimension ($d$). A thorough investigation of this relationship remains a future work.

\subsection{A candidate distribution}
Once we have experimentally established that strong universality holds,  a natural follow-up question is whether the observed limiting distribution $\cL^*$ is a familiar one, and in particular if it has a simple expression. Surprisingly, it seem that the answer might be yes.
We denote the \emph{left-skewed Gumbel distribution}  by $\mathrm{LGumbel}$, whose CDF and PDF are given by
\eqb\label{eqn:gumbel}
F(x) = 1-e^{-e^x},\quad\text{and}\quad f(x) = e^{x-e^{x}}.
\eqe
The expected value of this distribution is precisely the Euler-Mascheroni constant ($\lambda$) that we use in \eqref{eqn:AB}.
In Figures \ref{fig:loglog} and \ref{fig:noniid} the black dashed lines represent the CDF and PDF of the LGumbel distribution. In addition, the left column of these figures presents the QQ-plot of all the different models compared to the LGumbel distribution. These plots provide a very strong evidence for the validity of our next conjecture.

\begin{con}\label{con:gumbel}
$\cL^* = \mathrm{LGumbel}$.
\end{con}

Note that if $Z$ is an LGumbel random variable, then $e^Z$ is a standard exponential random variable (i.e., Exponential(1)). Therefore, if the distribution of $\set{\ell(p)}_{p\in \dgm_k}$ is approximately LGumbel, then the distribution of $\{(\log (\pi(p)))^A\}_{p\in\dgm_k}$ is approximately $\mathrm{Exponential}(\lambda_{\S,\cT,k})$ for some $\lambda_{\S,\cT,k}>0$ which  can be easily estimated from the data, but otherwise is unknown at the moment (in fact $\lambda = e^B$, where $B$ is explored in Figure \ref{fig:means}). We note that in the case where $\P_n$ is an iid measure, Conjecture \ref{con:iid} implies that $\lambda_{\S,\cT,k}$ depends on $d$ but is otherwise independent of $\S$.

\subsection{Independence}
In Conjectures \ref{con:iid} and \ref{con:loglog}, we discussed the limit for the empirical measure of  persistence $\pi$-values and $\ell$-values, respectively.
The way we interpret these conjectures is as follows. Given a random sample $\bX_n \sim \P_n$, the persistence diagram $\dgm_k$  is a random point process in $\R^2$, i.e.~a collection of random points $\set{p_1,\ldots, p_M}\subset \R^2$, where $M$ itself is random. Note that there is no natural ordering to the points in a diagram, and to some extent their ordering is algorithm-specific. Therefore, in our point process model we will assume that the points are ordered completely at random (i.e.~whatever the output of the algorithm is, we apply a random permutation to the points). In this case, one can show that given $M=m$,  the random variables $\pi(p_1),\ldots, \pi(p_m)$ all have the same marginal distribution. This distribution is approximated by the empirical measures we studied earlier (i.e., $\Pi_n$ for $\pi(p_i)$, and $\cL_n$ for $\ell(p_i)$). However, even if we knew exactly what the distribution of $\pi(p_i)$ is, it would not tell us anything about the \emph{joint distribution} of all $\pi$-values together, and in particular about the dependency relationship between them.

\begin{figure}[h]
\begin{center}
	\begin{subfigure}{0.45\linewidth}
	\centering
 		\includegraphics[width=0.7\linewidth]{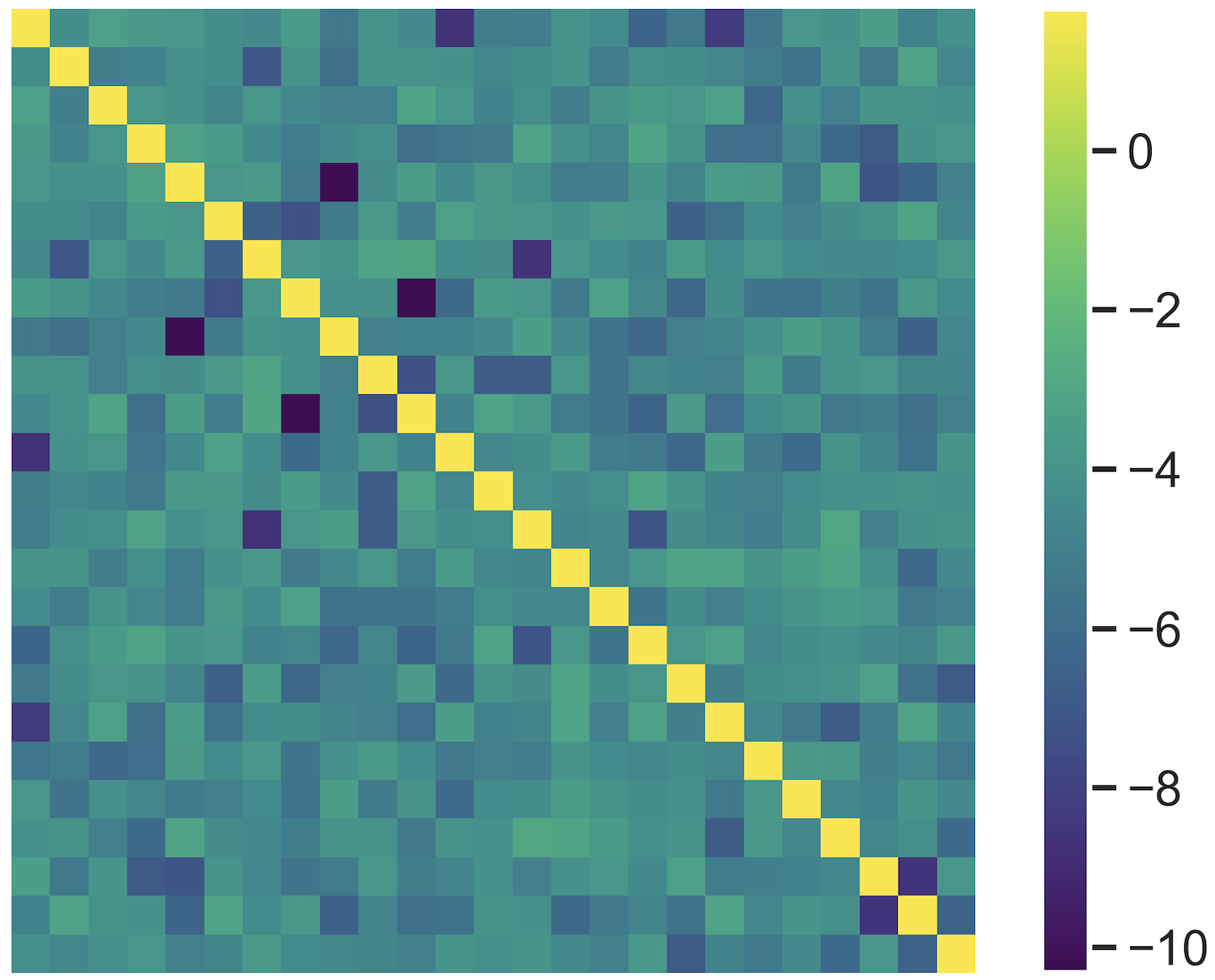}
\end{subfigure}
		\begin{subfigure}{0.45
		\linewidth}
	\centering
		\includegraphics[width=0.7\linewidth]{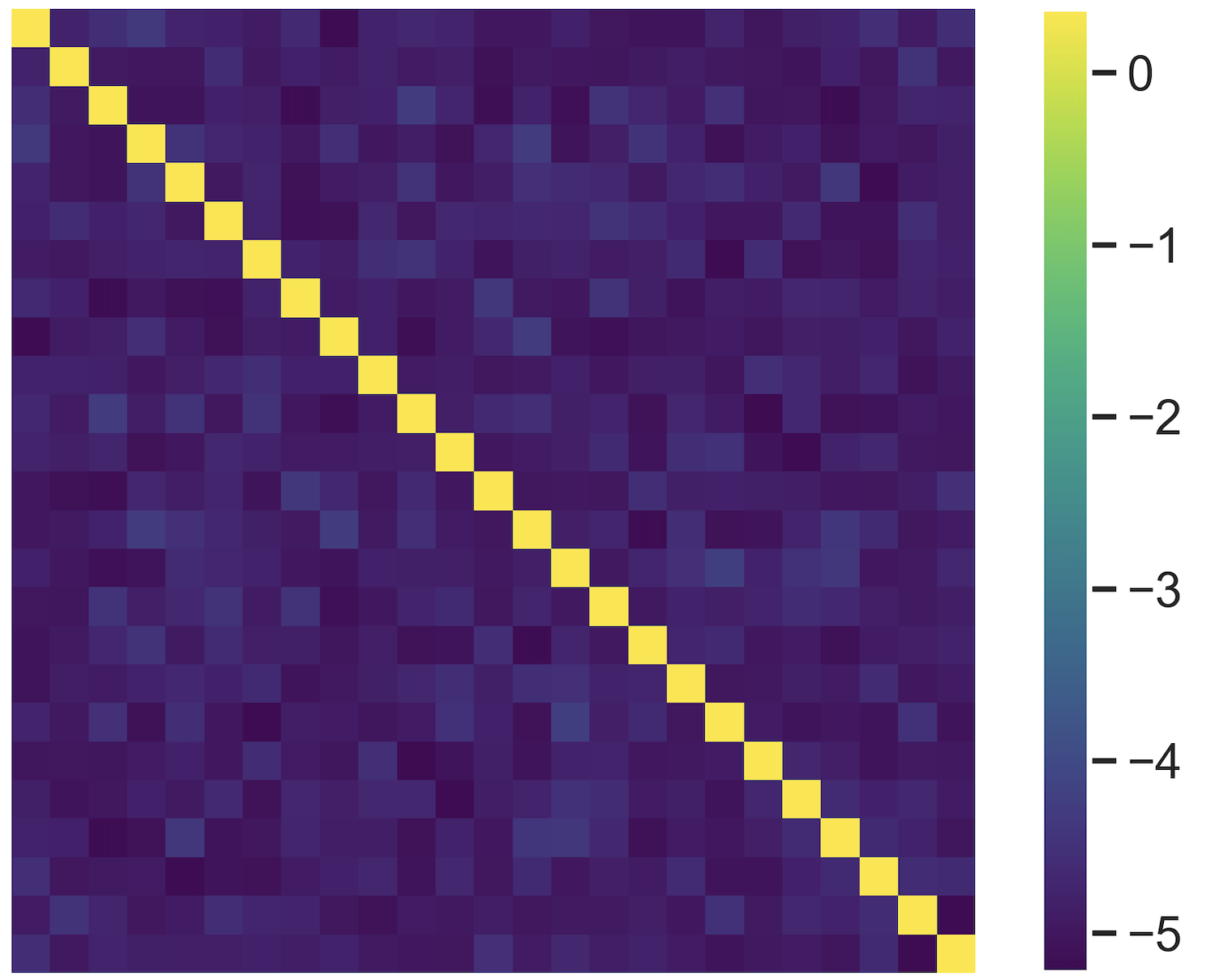}
	\end{subfigure}
	\end{center}
 	\caption{(left) The estimated covariance function for $\ell(p_1),\ldots,\ell(p_{25})$, generated by $N=100,000$ repeated experiments. (right)The estimated  distance-covariance function.  The resulting matrices are similar to the iid case.
 	\label{fig:indep}}
\end{figure}

In this section we study the dependency between cycles within a single diagram. For simplicity, we focus on two simple cases here -- a uniform iid sample in a 2-dimensional box and a sampled 2-dimensional Brownian motion. The results were quite similar and therefore we present only the 2-dimensional box here, while the Brownian motion results are available in Appendix \ref{app:indep}.
The first experiment was to estimate the covariance matrix for $\ell(p_1),\ldots, \ell(p_{25})$, estimated from $N=100,000$ independent persistence diagrams. The result is presented in Figure \ref{fig:indep} (left). The structure of the covariance matrix highly suggests that the points are uncorrelated. The average correlation coefficient, in absolute value, among all $300$ possible pairs of variables is 0.002, and the maximum (in absolute value) is 0.009. Not only are these values quite low, they are comparable to the values obtained in an empirical covariance matrix obtained by generating iid random variables from the LGumbel distribution (see Appendix \ref{app:indep}).

To further test the dependency among persistence values, we compute the matrix of distance covariance values. For a given pair of random variables, the distance covariance value \cite{szekely2007measuring} estimates how far
the joint characteristic function of the variables is from the product of the marginal characteristic functions. The resulting matrix is presented in Figure \ref{fig:indep} (right). The average distance correlation is 0.005, and its maximum is 0.01. These values are also comparable to iid random variables.

In conclusion, while our experiments above only consider pairwise dependency, they provide a strong indication that the $\ell$-values are independent (in some asymptotic sense). For more details, see Appendix \ref{app:indep}.

\section{Experiments}\label{sec:exp}
In this section we provide more details about the large body of experimental evidence we collected to support 
the conjectures stated in the previous section.
We wish present the highlights of our results, in a brief yet clear way. For elaborate details, as well as complete set of results, we refer the reader to  Appendix \ref{app:exp}.
The experiments were run across a number of dimensions (typically $d=2,3,4,5$), and across a range of sample sizes ($n\in[10^4,10^6]$). 
As the results  do not seem to be significantly impacted by the choice of $n$ (within the above range), we leave this detail to the appendix.

\begin{remark}
While we phrase all the results here in terms of the \v Cech complex, in practice the computations were done using the much lighter Delaunay-Alpha complex, which has an identical persistence diagram. 
\end{remark}
 \begin{remark}
The experiments were run primarily on the Apocrita HPC cluster, using several different software packages including Gudhi~\cite{maria2014gudhi}, Ripser~\cite{bauer2021ripser}, Eirene~\cite{henselmanghrist16},  Dionysus~\cite{dionysus}, and Diode~\cite{diode}. We plan to publish the code as well as the datasets used in this paper. The details will appear in the next revision of this paper.
\end{remark}

\subsection{Sampling from iid distributions}
The most natural samples to start with are those where the points $X_1,\ldots,X_n$ are independent and identically distributed (iid). Here, we can test  Conjecture \ref{con:iid} as well as Conjectures \ref{con:loglog} and \ref{con:gumbel}. 

We started by considering samples from the uniform distribution on various compact manifolds (with and without a boundary), with diverse geometry and topology. See Table \ref{tab:iid_uniform} for the list of settings tested. Note that for the torus, Henneberg surface, Klein bottle and projective plane, our sample is not uniformly distributed in the manifold itself, but rather in the space of its given parametrization.
The results for most of the settings are included in Figures  \ref{fig:pi_cdfs} and \ref{fig:loglog}. In order to avoid overloading these figures, we omitted some of the cases, while the results for all settings can be found in  Appendix \ref{app:iid}.

Next, we wanted to test non-uniform distributions (in $\R^d$).
To this end we took three representative distributions in $\R$, and for each distribution $\cD$ we generated random points $X_i = (X_i^{(1)},\ldots, X_i^{(d)})$ where the coordinates are independent, and each $X_i^{(j)}$ is distributed according to $\cD$. The  beta distribution, is non-uniform in a compact space ($[0,1]$). The normal and Cauchy distributions are supported on $\R$. 
The standard normal distribution is an obvious choice and Cauchy was chosen as a heavy-tailed distribution without any moments. The details are in Table \ref{tab:iid_nonuniform}. 

All the settings mentioned above have a well-defined dimension $d$, and indeed support Conjecture \ref{con:iid} (for different dimensions $d$), as well as Conjectures \ref{con:loglog} and \ref{con:gumbel}.

In addition to the standard distributions considered above, we also tested two more complex models (still iid). {\bf Stratified spaces:} These are mixed-dimensional spaces, as suggested in \cite{thomas2021learning}. As an example, we  considered a $2$-dimensional square embedded in a $3$-dimensional cube. For each $i$, with probability $p\in(0,1)$ we sampled $X_i$ uniformly from the square, and with probability $(1-p)$, we sampled it uniformly from the cube.
{\bf Linkages:} We sample randomly from the 
configuration space of a closed five-linkage with unit length, i.e. regular pentagons with unit side length (up to translations and rotations). This configuration space is a 2-dimensional subspace of $\mathbb{R}^6$.
The results for both models seem to support Conjectures \ref{con:loglog} and \ref{con:gumbel} (see Appendix \ref{app:iid}). 
We note that for the stratified space the value of $B$  (\ref{eqn:loglog}) varies with $p$.

\begin{table}[]
    \centering
    \begin{tabular}{|c|c|c|}
    \hline
        \multirow{2}{*}{space ($\cS$)}&
        intrinsic& 
        embedding 
        \\
        &dimension ($d$)& dimension \\
        \hline 
        \hline
             Box & 2,3,4,5& $d$ \\
                Ball & 2,3,4,5& $d$\\
                Sphere& 2,3,4 & $d+1$\\
                Torus & 2 & 3  \\
                Henneberg surface & 2 & 3 \\
                Neptune mesh & 2 & 3 \\
                Klein bottle & 2 & 4 \\
                Projective plane & 2 & 4 \\ \hline
    \end{tabular}
    \caption{An exhaustive list of all iid uniform distributions tested. In the \v Cech complex, we computed all relevant dimensions of homology, $k=1,\ldots, d-1$. For computational reasons, we limited the Rips to $k=1,2$.
    \label{tab:iid_uniform}}
\end{table}

\begin{table}[]
    \centering
    \begin{tabular}{|c|c|c|c|}
    \hline
        distribution ($\P$)&
        density function &
        space ($\cS$)&
        dimension ($d$)\\
        \hline
        \hline
            Beta(3,1) & $x^2$ & $[0,1]^d$ & 2,3,4,5\\
             Beta(5,2) & $x(1-x)^4$ & $[0,1]^d$ & 2,3,4,5 \\
                Normal & $\exp(-x^2/2)$ & $\R^d$ & 2,3,4,5\\
                Cauchy & $(1+x^2)^{-1}$ & $\R^d$ & 2,3,4,5  \\
                
                         \hline
    \end{tabular}
    \caption{An exhaustive list of all iid non-uniform distributions tested. The homology dimensions computed ($k$) are the same as in Table \ref{tab:iid_uniform}. The second column contains the expression of the 1-dimensional density used  (up to a normalizing constant).
    \label{tab:iid_nonuniform}}
\end{table}

\subsection{Sampling from non-iid distributions}
The iid settings are the most commonly studied  in the stochastic topology literature, and are therefore the first models to examine. In order to gain a better understanding about the extent of universality, as well as to consider more realistic models, we wish to test more complex dynamics as well.

We tested two vastly different models. {\bf Brownian motion:} The Brownian motion is a real-valued continuous-time Gaussian process with stationary independent increments. For $t\ge 0$ we define the $d$-dimensional Brownian motion  $W_t = (W_t^{(1)},\ldots, W_t^{(d)})$ where the coordinates $W_t^{(i)}$ are independent real-valued Brownian motions. The point-cloud we generate here is a discrete-time sample of $W_t$, at times $t=1,\ldots,n$. In other words, for $1\le i\le n$ we define $X_i = W_{t=i}$.
Note that the variables $X_1,\ldots,X_n$ are neither independent, nor identically distributed. In addition, note the path of Brownian motion path has a fractal dimension. 
{\bf The Lorenz system:}
Here we took a discrete-time sample of the Lorenz dynamical system -- a well-studied  example for a chaotic  system.
It is a 3-dimensional system of non-linear differential equations, given by:  $\dot{x} = \sigma (y-x), \dot{y} = x(\rho - z) - y, \dot{z}=xy-\beta z$, for some positive $\sigma,\rho,\beta$.

In Figure \ref{fig:noniid} we present the Brownian motion and Lorenz system examples. More cases can be found in Appendix \ref{app:non-iid}.
Quite surprisingly, the persistence diagrams generated by  these complex systems seem to follow Conjectures \ref{con:loglog} and \ref{con:gumbel} as well.

\subsection{Testing on real data}
The ultimate test for Conjectures \ref{con:loglog} and \ref{con:gumbel} is in real world data. We tested two different examples.
{\bf Natural images:} 
The image database is taken from van Hateren and van der Schaaf \cite{van1998independent}. This dataset is a large collection of natural gray-scale images.
From these images, we randomly extract patches of size $3\times 3$, forming a 9-dimensional point-cloud. We applied the dimension reduction procedure proposed by \cite{lee_nonlinear_2003}, which results in a point-cloud on a $7$-dimensional sphere embedded in $\R^8$. We tested both the $7$-dimensional point-cloud, as well as its lower-dimensional projections.
{\bf Audio recording:} 
We took an arbitrary speech recording, and applied the time-delay embedding transformation (cf.\cite{takens_detecting_1981}) to convert the temporal signal into a $d$-dimensional point-cloud. More details are in Appendix \ref{app:real-data}.

We tested both point-clouds in various dimensions, using both the \v Cech and Rips complexes. The matching to the universal distribution was quite remarkable. The results are in Figure \ref{fig:noniid}.

\begin{remark}
While nearly every point-cloud we tested supported Conjectures \ref{con:loglog} and \ref{con:gumbel}, we found two exceptions. One is when we take points on a grid with small random perturbations, and the other is the Ginibre ensemble, which is an example of a repulsive determinantal point process. In both cases the spacings between points are quite homogeneous, which makes it impossible for clusters to form. This has a significant effect on the distribution of large $\pi$-values, since  those are generated by dense clusters. 
\end{remark}

\begin{remark}
In all the figures and computations we present in this paper, in theory we have to assume that the diagrams we are looking at contain only noisy cycles and no signal. This is important, for example, if we want to estimate the CDF, or the value of $B$ \eqref{eqn:AB}. In most cases, where the signal is known (i.e., an annulus) we can easily remove these points manually. However, in cases where it is not known in advance (i.e., audio sampling) we cannot do so.
However, in practice, this is not a significant issue. The diagrams we compute contain between thousands to millions of points. Thus, as long we do not expect the signal to contain more than a few points (as in most realistic examples), even if we include the signal cycles in our estimates, their effect will be negligible, especially in the $\logg$ scale. This is especially important for the applicability of our hypothesis testing framework.
\end{remark}

\section{Applications}
In this section we present the most natural use for our Conjectures in the context of hypothesis testing for persistence diagrams. We divide the discussion between finite and infinite cycles, and present a few examples at the end.
\subsection{Hypothesis testing}\label{sec:ht_fin}
Suppose that we are given a persistence diagram $\dgm = \set{p_1,\ldots,p_m}$. Our goal is to determine for each point $p_i$ whether it is signal or noise. This can be modelled a multiple hypothesis testing problem with the $i$-th null-hypothesis we are testing, denoted $H_0^{(i)}$, is that $p_i$ is a noisy cycle. Assuming Conjectures \ref{con:loglog} and \ref{con:gumbel} are true, we can formalize the null hypothesis in terms of the $\ell$-values (\ref{eqn:loglog}) as follows
\[
H_0^{(i)}\ :\ \ell(p_i) \sim \mathrm{LGumbel}.
\]
In other words, cycles that deviate significantly from the LGumbel distribution should be declared as \emph{signal}. Given $\ell(p_i) = x$ the $i$-th p-value, using the LGumbel distribution, is 
\eqb\label{eqn:pval}
\text{p-value}_i = \cprob{\ell(p_i) \ge x}{H_0^{(i)}} = e^{-e^x}.
\eqe

Since we are testing multiple cycles simultaneously, the p-values should undergo some correction first. 
The simplest to use is the Bonferroni correction, which seems to suffice in our experiments. Alternatively, more accurate corrections may be used, especially if the independence conjecture we discussed earlier holds. To conclude, the signal part of a diagram could be recovered with signficance level $\alpha$ via
$$ \dgm_k^{_\mS}(\alpha) = \set{p\in \dgm_k :  e^{-e^{\ell(p)}} < \frac{\alpha}{|\dgm_k|}}.$$

\subsection{Infinite cycles}\label{sec:ht_inf}
Computing persistent homology for either the \v Cech or Rips filtrations, it is rarely the case that one  computes the entire persistence diagram, i.e., taking the entire range of radii. Since large radii introduce a large number of simplices, computing persistent homology becomes intractable.
Therefore, it is required to heuristically choose  a threshold radius $\tau$ for the persistence computation. We denote the resulting diagram $\dgm_k(\tau)$. In this case,   one often discovers cycles that are born prior to $\tau$, but die after $\tau$. From the algorithm perspective  such cycles are born but never die, and hence are declared ``infinite''.
The question we try to address in this section is how to \emph{efficiently} determine whether such infinite cycles are  statistically significant or not. 

Suppose that $p=(\mathrm{b},\mathrm{d})\in\dgm_k(\tau)$ is such an infinite cycle, i.e., $\mathrm{b}\le\tau$ is known and $\mathrm{d}>\tau$ is unknown.
While we cannot compute $\ell(p)$, we can still determine that $\ell(p)> \tau/\mathrm{b}$. This, in turn, provides us with an upper bound for the p-value of $p$
\[
    \text{p-value}_i  < \text{p-value}_i(\tau) := e^{-e^{\tau/\mathrm{b}}}.
\]
There are now two possible options. If $\text{p-value}_i(\tau)$ is below the required significance value (e.g.~$\alpha/|\dgm_k|$), we can determine that $p$ is indeed significant, without knowing its true death time. Otherwise, we can determine the minimal value $\tau^*$ required so that $\text{p-value}_i(\tau^*)$ is below the significance value. We can then compute persistent homology again, with $\tau^*$ as the new threshold. If the cycle represented by $p$ remains infinite (i.e. $\mathrm{d} > \tau^*$) we now declare it as significant. The key point here is that the actual death time $\mathrm{d}$ can be much larger than $\tau^*$. However, for the sole purpose measuring significance, we do not need to know the exact value of $\mathrm{d}$, only whether it is smaller or larger than $\tau^*$.

The procedure we just describe works well, if we are only studying a single infinite cycle. However, it is very likely that we find more than a single infinite cycle in $\dgm_k(\tau)$. Moreover, once we increase the threshold from $\tau$ to $\tau^*$ it is possible that new infinite cycles will emerge as well. We therefore propose the iterative procedure described in Algorithm \ref{alg:cap}. Briefly, this algorithm looks at every step for all infinite cycles, and chooses the next threshold $\tau$ so that we can determine whether the earliest-born infinite cycle is significant or not. Note that it is possible instead to  take the latest-born infinite cycle (switching $\min(I)$ with $\max(I)$). This will make sure that \emph{all} currently infinite cycles are determine in the next iteration, and thus will reduce the total number of iterations. However, it is also possible that the new threshold $\tau$ computed this way will be an overshoot, and will lead us to generate a filtration much larger than we actually need. In principle, one can choose any of the values in $I$ in order to determine the next threshold, in a way that balances between the number of iterations   and the size of the filtration used.
 
The value $\pi_{\min}(x)$ in the Algorithm \ref{alg:cap}, is the minimum $\pi$-value required so that the resulting p-value (\ref{eqn:pval})  is smaller than $x$. Formally,
\[
\pi_{\min}(x) = \ell^{-1}\left(F^{-1} \left( 1-x\right)\right) = \ell^{-1}(\logg (1/x)),
\]
where $F$ is the CDF of the LGumbel distribution (\ref{eqn:gumbel}).

\begin{algorithm}[h!]
\caption{Finding threshold for infinite cycles}\label{alg:cap}
\begin{algorithmic}
\State $\tau \gets \tau_0$
\Do
    \State $\cD \gets \dgm_k(\tau)$
    \State $I \gets \set{\mathrm{b}: (\mathrm{b},\mathrm{d})\in \cD, \mathrm{d}=\infty, \text{ and } \tau/\mathrm{b} < \pi_{\min}(\alpha/|\cD|)}$
   \State $\tau \gets \begin{cases} \min(I)\cdot \pi_{\min}(\alpha / |\cD|) & I \neq \emptyset\\ \tau & I = \emptyset \end{cases}$
\doWhile{$|I|>0$}
\State \Return $\tau$
\end{algorithmic}
\end{algorithm}

\begin{remark}
Note that since the total number of cycles is non-decreasing with $\tau$, we are not guaranteed to reduce the size of $I$ at every step. However, since $\tau$ is always increasing, the algorithm must eventually terminate for finite point-clouds. In the future we plan to look at refinements of this algorithm.
\end{remark}

\begin{remark}
Using the partial diagram $\dgm_k(\tau)$, an important caveat is that $|\dgm_k(\tau)| < |\dgm_k|$, which  affects the significance correction we use. The partial diagram may also lead to an erroneous computation of the $B$ value in \eqref{eqn:loglog}.
To minimize the effect of these issues, it is  imperative  that most of the noisy cycles are included in $\dgm_k(\tau)$. The iterative nature of Algorithm \ref{alg:cap} is aimed to achieve this purpose. In fact, the only way in which the algorithm will terminate without including all noisy cycles is if it happens to select a value of $\tau$ that does not cover all noisy cycles and at the same time does not introduce any infinite cycles. The probability, however, to select such $\tau$ is very small.
\end{remark}

\subsection{Examples}
We present a few examples for the hypothesis testing framework we presented in Sections \ref{sec:ht_fin} and \ref{sec:ht_inf}. In all our experiments, we set the desired significance level to be $\alpha=0.05$.

\begin{figure*}[h]
\begin{center}
	\includegraphics[width=\textwidth]{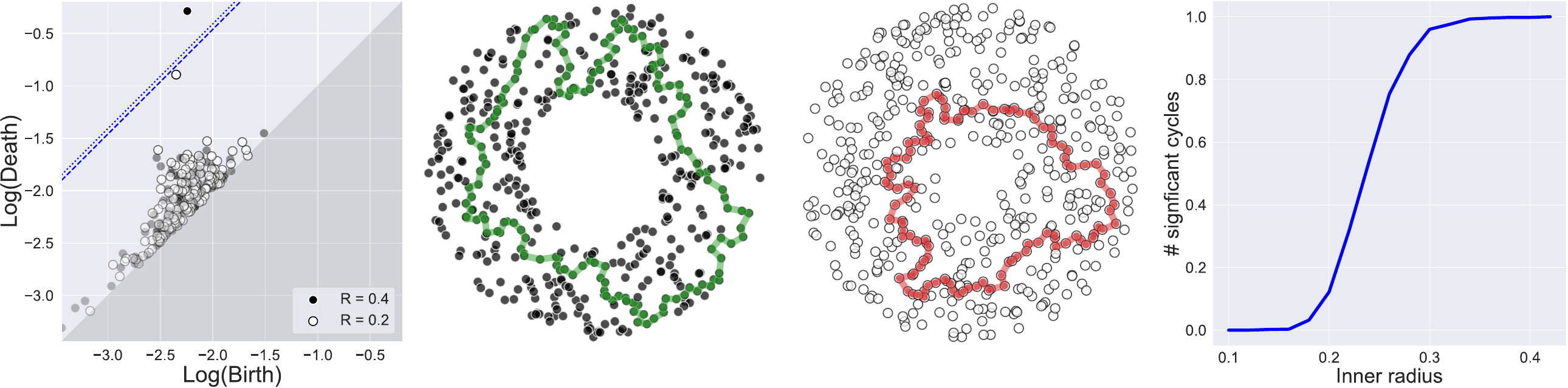}
	\end{center}
 	\caption{Computing p-values for an annulus. (left) Persistence diagrams for two instances with different inner radii, with the significance lines  shown, for $\alpha=0.05$. The dotted line shows the significance threshold for the diagram coming from the larger inner radius, and the dashed line is the threshold for the smaller radius. (middle-left) A significant cycle (p-value = 0.008) is detected for $R=0.4$. (middle-right) 
 	The most persistent cycle for $R=0.2$ is not detected as signal (p-value = 0.076), due to the small inner radius.
 	(right) The p-value curve as a function of the inner radius.\label{fig:annulus_full}}
\end{figure*}

\noindent {\bf Computing p-values.} We start with a toy example to test  our p-value computation. We sampled 1,000 points on an annulus in $\R^2$, whose outer radius is $1$ and the inner radius is varied. In this case we expect to detect a single $1$-cycle, corresponding to the hole of the annulus. For each value of $r$ we computed the persistence diagram, and counted the number of cycles that are declared as significant (i.e., $\text{p-value}<\frac{\alpha}{|\dgm|}$). 
In Figure~\ref{fig:annulus_full} on the right we present the results of this experiment  with the curve showing the average number of significant cycles detected, over 1,000 repetitions, versus the inner radius of the annulus. 
In this figure we also show the persistence diagrams (left) of two realizations with different radii -- one where the cycle is detected as significant, and one where it is not. For each diagram, the line of significance (corresponding to $\alpha=0.05$) is shown (the dotted and dashed lines correspond to the larger and smaller radii respectively).
The figures in the middle show the cycles with the largest $\ell$-value  for $R=0.4$ (left, significant) and $R=0.2$ (right, non-significant).  

\begin{figure*}[h]
\begin{center}
	\includegraphics[width=\textwidth]{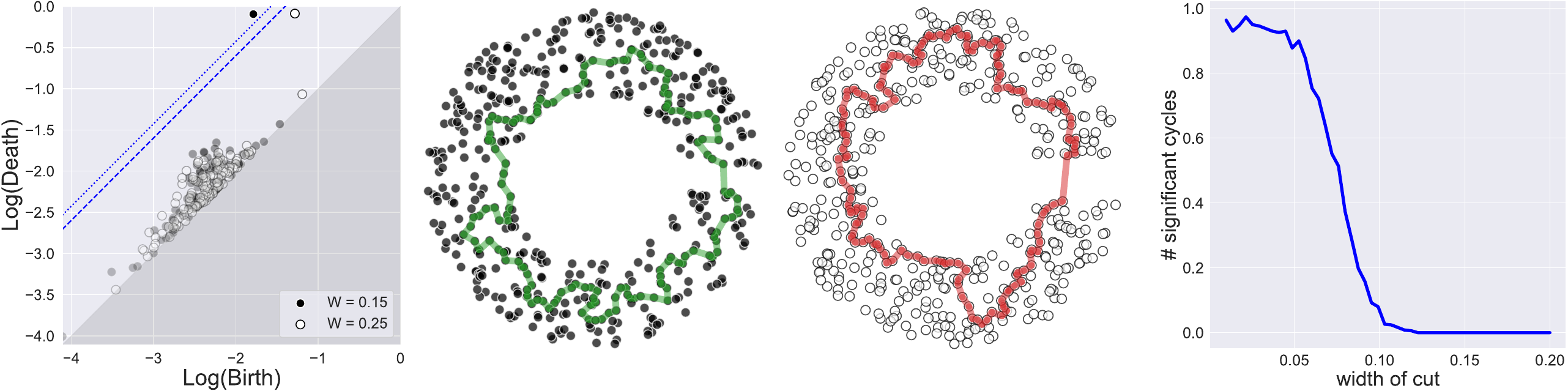}
	\end{center}
	
 	\caption{Computing p-values for cut-annulus. (left) Persistence diagrams for two instances with different cut-widths, with the significance lines shown for $\alpha=0.05$ (the dotted and dashed lines correspond to  the smaller and larger cuts respectively).
 	(middle-left) At cut-width $W=0.15$  it is difficult to distinguish the cut-annulus from a complete one, and the cycle is detected as signal (p-value = 0.008). 
 	(middle-right) A non-significant but persistent cycle (p-value = 0.16), with a cut-width of $W=0.25$. The lack of significance is due to the later birth time of the cycle. (right) The p-value curve for a cut-annulus with a function of the cut-width. \label{fig:pval_cut}}
\end{figure*}

\begin{figure*}[h]
\begin{center}
		\includegraphics[width=\textwidth]{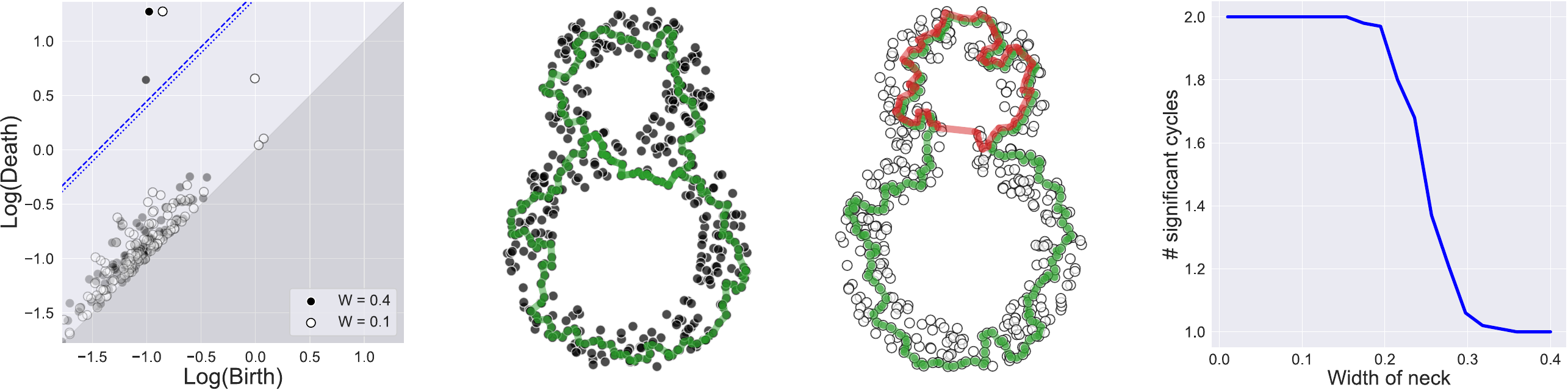}
	\end{center}
 	\caption{Computing p-values for a figure-8. (left) Persistence diagrams for two instances of the figure-8 with different neck gaps, with the significance lines shown for $\alpha=0.05$ (the dotted and dashed lines correspond to  the narrow and wider neck gap respectively) (middle-left) Two significant cycles (p-values = 0.0005,0.012) are detected with smaller neck gap ($W=0.1$). (middle-right) For $W=0.4$ the diagram contains two persistent cycles, but only the green one is declared as significant (p-value=0.0015), while the red one is not (p-value=1).
 	The lack of significance is due to the later birth time of the cycle. (right) The p-value curve for a figure-8 as a function of the neck gap.\label{fig:pval_fig8}}
\end{figure*}

In the next experiment, we used a ``cut'' annulus, i.e., we take an annulus with a strip removed from it (see Figure \ref{fig:pval_cut}). The resulting space is contractible, and in particular has no signal to be detected. However,  the point-cloud generated on the cut-annulus is very similar to that generated on a complete annulus. Indeed, the persistence diagram for points on the cut-annulus tend to have a single 1-cycle with a long lifetime
(see Figure \ref{fig:pval_cut}, left). Our goal is to show that, using our p-value computation,  we can determine that despite  its long lifetime, this cycle should not be detected as signal. 
In Figure~\ref{fig:pval_cut}  (middle left and middle right) we present the largest cycles for two different cut-widths.
While these cycles look quite similar,
the one on the right is not considered statistically significant, due to its later birth time.  In Figure \ref{fig:pval_cut}(right), we show the average number of signal cycles detected in the cut-annulus as a function of the width of the cut. When the width is very small (as in the middle-left figure), it is impossible to distinguish between the annulus and cut-annulus with just 1,000 points. Indeed, in these cases, we often we falsely detect a signal. However, already for relatively small values (around $W=0.05$) there are cases where we manage to reject the ``fake'' signal, and at around $W=0.15$ we do so perfectly. In other words, using our p-values we can reject cycles, even with relatively long lifetimes, that would otherwise be falsely detected as signal. Finally, in Figure~\ref{fig:pval_fig8} we show a similar figure-8 experiment, where we vary the width of the neck.

\begin{figure*}[h]
\begin{center}
	\centering
 		\includegraphics[width=0.8\textwidth]{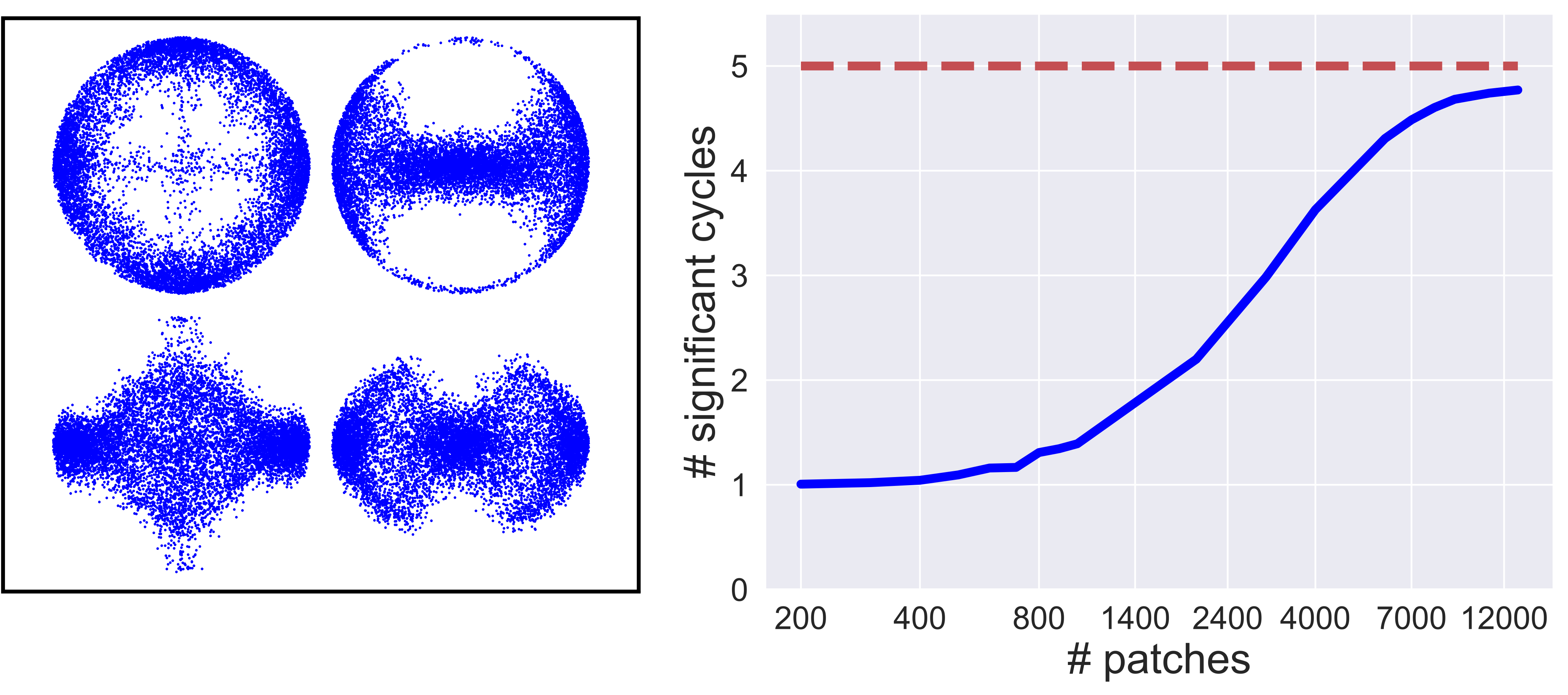}
	\end{center}
 	\caption{ Testing the 3-circle model in the natural image patches. (left) Four different 2-dimensional projections of the 8-dimensional patches point-cloud. The projection on the top-right shows a 1-cycle that is quite ``thin'' and contains large gaps. (right) The p-value curve for the patches dataset as a function of the number of samples. We observe that we are not detecting 5-cycles in 100\% of the cases. This is most likely due to the ``thin'' cycle.
 	\label{fig:patchespvals}}
\end{figure*}

Finally, we apply this to a real-world dataset, specifically the van-Hateren natural image database mentioned earlier. The main claim in \cite{de_silva_topological_2004} is that the space of $3\times 3$ patches (after dimension reduction, normalization, and filtering), has a 3-circle structure in $\R^8$. This claim was later used to deduce that the space of patches is concentrated around a Klein-bottle \cite{carlsson_local_2008}. The main support for this claim, is the persistence diagram computed over the patches, where five  1-cycles have a relatively long lifetime. Here we use our p-value computation to provide a quantitative statistical support for this claim.
To this end, we randomly selected a subset of the patches, processed it as in \cite{de_silva_topological_2004}, and then computed p-value for the cycles in the persistence diagram (using the Rips complex). We repeated this experiment for varying number of patches, and computed the average number of detected signal cycles over 250 trials. The results are presented in Figure~\ref{fig:patchespvals}. Firstly, we observe that there exists a single 1-cycle that is always detected, making it a very strong signal. As we increase the sample size, the number of detected signal cycles increases as well, and indeed gets closer to five.
We observe, however, that one of the five cycles is not always detected. when we plot a 2-dimensional projection of the points, we see indeed that this cycle is made of points with 
a low density, that contain relatively large gaps. Such large gaps increase the birth time, and consequently the p-value as well. To conclude, using our hypothesis-testing method we are not only able to correctly detect the signal cycles discovered in \cite{de_silva_topological_2004}, but we can also quantitatively declare how significant each of these findings is.

\noindent {\bf Infinite cycles.}
Here, we wanted to demonstrate the power of Algorithm \ref{alg:cap}.
We generated a point-cloud on the 2-dimensional torus, embedded in $\R^3$. The torus radii are $R_1=1$ and $R_2=2$. Originally, if we wanted to inspect the lifetime of the signal cycles in this example, and possibly to compute a p-value, we need to take the filtration all the way up to $r=R_2$. This will result in a large amount of simplices that would require a massive amount of memory and CPU time. Instead, we ran Algorithm \ref{alg:cap}, which incrementally increases the filtration threshold just enough to be able to determine significance for all the cycles in the filtration. Table \ref{tab:speedup} provides some insight into how much computational resources we can save by using this method. For varying sample sizes, we ran this procedure and found the largest value of $\tau$ needed to classify all cycles as signal or noise.  We compare the number of edges in the Rips complex at radius $\tau$ with the number of edges at radius $2$ (which, without Algorithm \ref{alg:cap}, is the only way to compute a p-value for both signal cycles). We observe that the amount of edges we are saving is quite large, and is increasing with the sample size. The saving we expect to achieve in the number of triangles is much higher, however, we were unable to compute the exact numbers due to computational limitations.

\begin{table}[]
    \centering
    \begin{tabular}{|c|c|c|c|c|}
    \hline
        \multirow{2}{*}{ \# of points}&\multirow{2}{*}{$\tau$}& \# of edges at $r=2$ & \# of edges at $r=\tau$& \multirow{2}{*}{ ratio }\\
        &&(1000s) &(1000s)&\\
            \hline
               \hline
  2000	& 1.98&	490 &	475 &	1.03\\
5000	& 1.67&	3067 &	1876 &	1.64\\
10000	&1.41&	12221 &	5022 &	2.43\\
20000	&1.19&	49074 &	13942& 	3.52\\
50000	&0.95&	306580 & 53581& 	5.72\\
           \hline
    \end{tabular}
    \caption{The number of edges (in 1000s) at the full threshold ($r=2$) versus the computed threshold $\tau$ versus the number of points. The last column  represents the  ratio of the number of edges.}
    \label{tab:speedup}
\end{table}

\section{Discussion}
Persistent homology is a powerful topological tool  used for anaylzing the structure of data.   Revealing the distribution of persistence diagrams has been  one of the greatest challenges faced by probabilists and statisticians working in the field. In this paper we argue that  viewed through the lens of the $\ell$-values, almost all persistence diagrams obey the same universal distribution.
Te support this statement we provided an exhaustive set of experimental results, ranging from artificial iid samples to real data. We expect the outcomes of this paper to be twofold. Firstly, the conjectures in this paper will provide fertile grounds for whole new line of theoretical research in stochastic topology.
We expect that revealing the full extent of Conjectures \ref{con:iid}-\ref{con:gumbel} will be a long lasting process, that will require the development of novel techniques and methods.
Secondly, we expect  new powerful statistical tools to be developed based on the universality properties presented here. The hypothesis testing framework we suggested here is already quite powerful, but we expect it to be only the tip of the iceberg. 
The conjectures presented here, open the door to developing  novel statistical tests for automatically detecting structure in data. In addition, they can lead to stronger stability statements and tighter bounds on distances between persistence diagrams in cases where there is sufficient randomness. 
Our results  may also lead to further algorithmic advances in computing persistent homology. 
For example, taking advantage of universality can lead to optimizing the software heuristics used in a principled manner.
To conclude, proving the conjectures in this paper will provide an explicit expression for the noise distribution in persistence diagrams with very few model assumptions. This, in turn, will open up new horizons in TDA, and might provide solutions to long standing open problems in the area.

\vspace{5pt}
\noindent{\bf Acknowledgements. }
We would like to thank Henry Adams for useful advice and suggestions regarding some of the examples presented here. We would also like to thank the following people for useful discussions: Robert Adler, Matthew Kahle, Anthea Monod, and Sayan Mukherjee.

\bibliography{zotero}
\bibliographystyle{plain}

\appendix
\newpage
\section{Experimental Details and Complete Results}\label{app:exp}
In this section we wish to expand the discussion in Section \ref{sec:exp}, and to provide more details about the experiments, as well more results that were not included in the paper body.

\subsection{Sampling from iid distributions}\label{app:iid}
\begin{figure}[h!]
     \centering
     \includegraphics[width=0.95\textwidth]{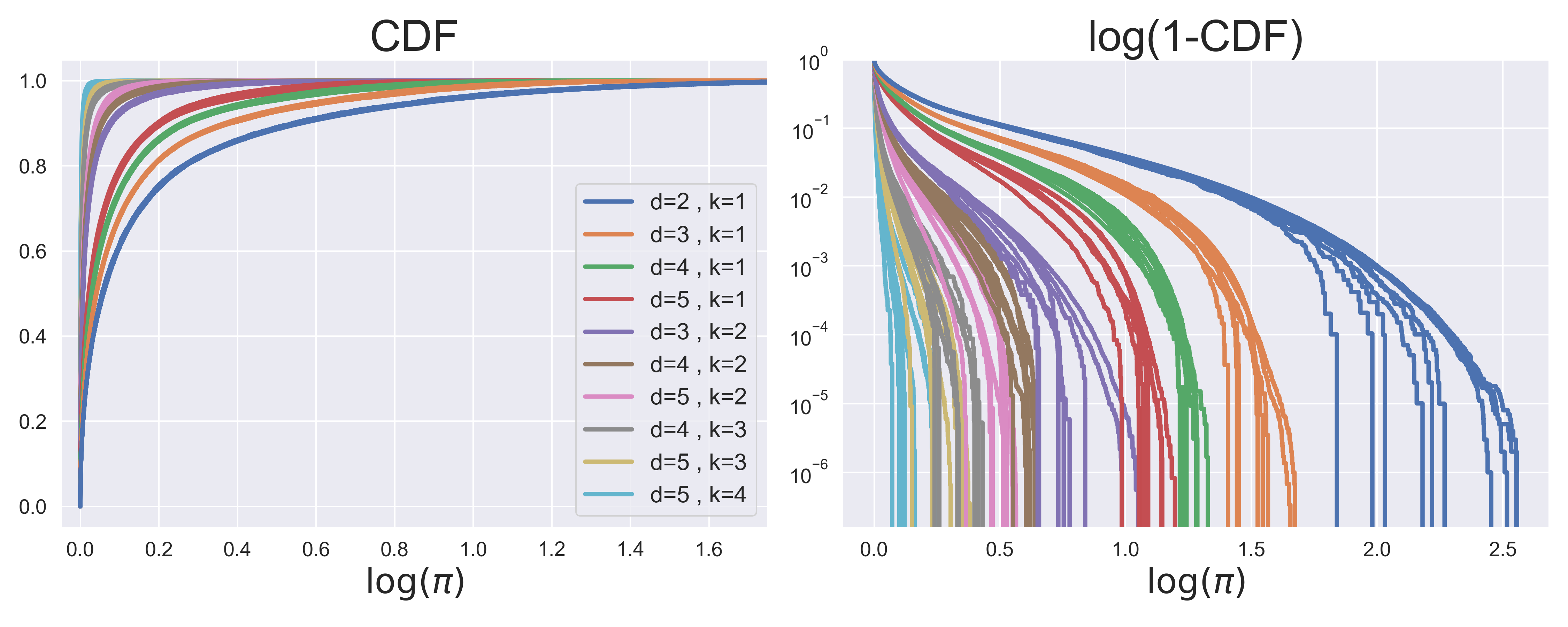}
     \caption{CDFs $\pi$-values in all iid settings tested, using the \v Cech complex. The complete list of models included in this figure is in Table \ref{tab:cech_iid_all}
     \label{fig:cech_iid}}
 \end{figure}

\begin{figure}[h!]
     \centering
     \includegraphics[width=0.95\textwidth]{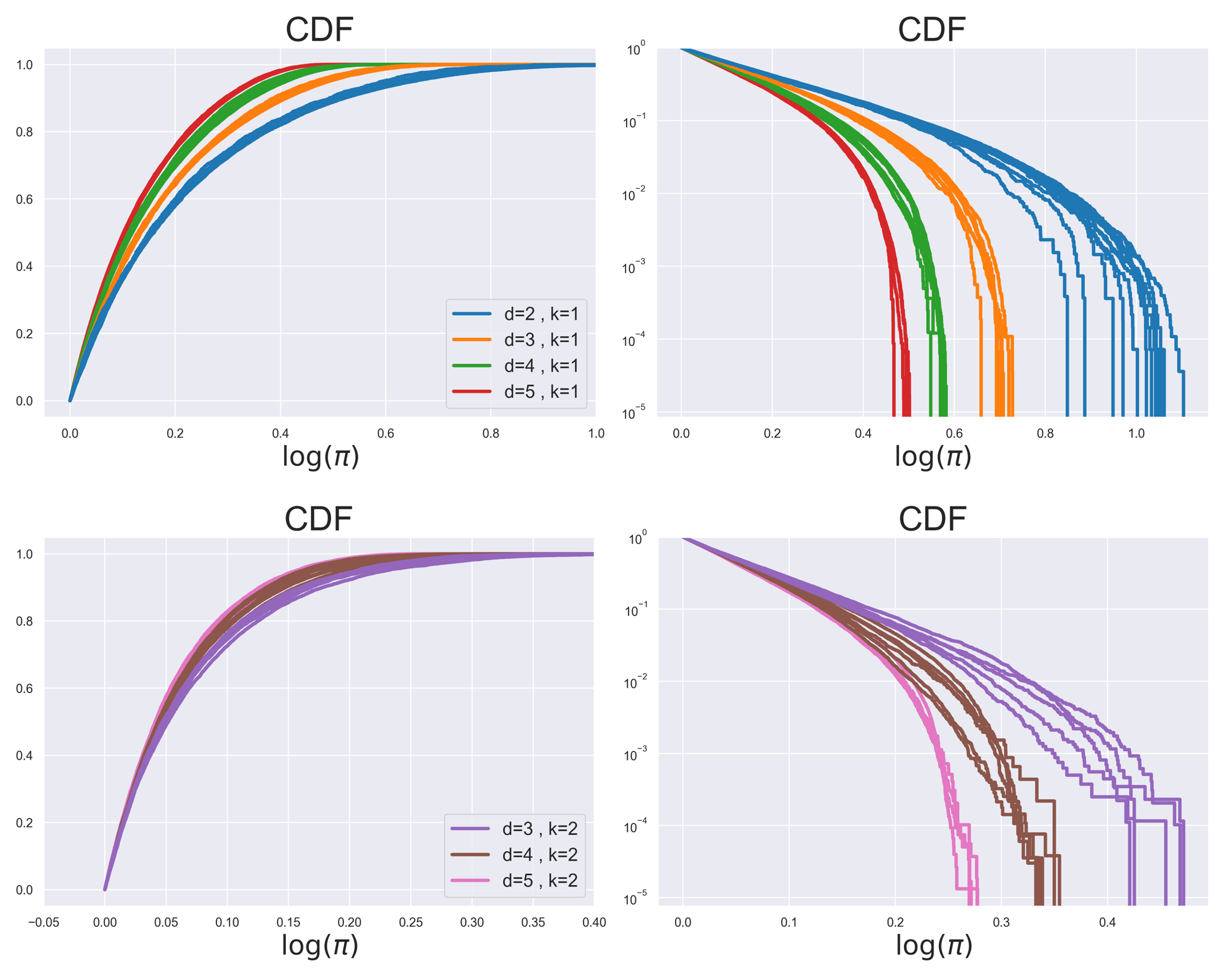}    
\caption{CDFs $\pi$-values in all iid settings tested, using the Rips complex. The complete list of models included in this figure is in Table \ref{tab:rips_iid_all}}
     \label{fig:rips_iid}
 \end{figure}

\begin{figure}[h!]
     \centering
     \includegraphics[width=0.95\textwidth]{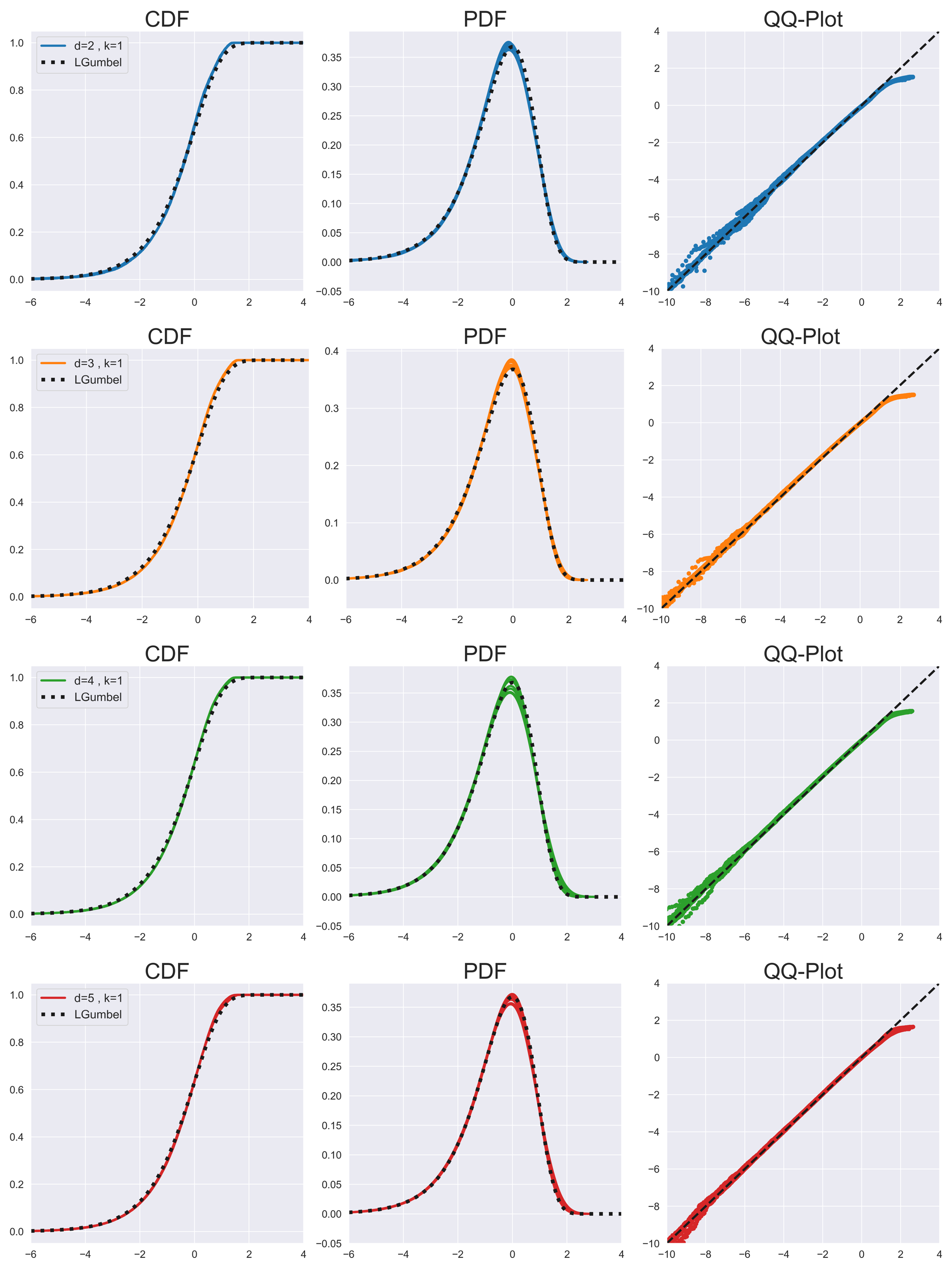}
      \caption{The $\ell$-values distribution for the \v Cech complex -- part 1/3.
    \label{fig:cech_iid_loglog_1}}
 \end{figure}

\begin{figure}[h!]
     \centering
      \includegraphics[width=0.95\textwidth]{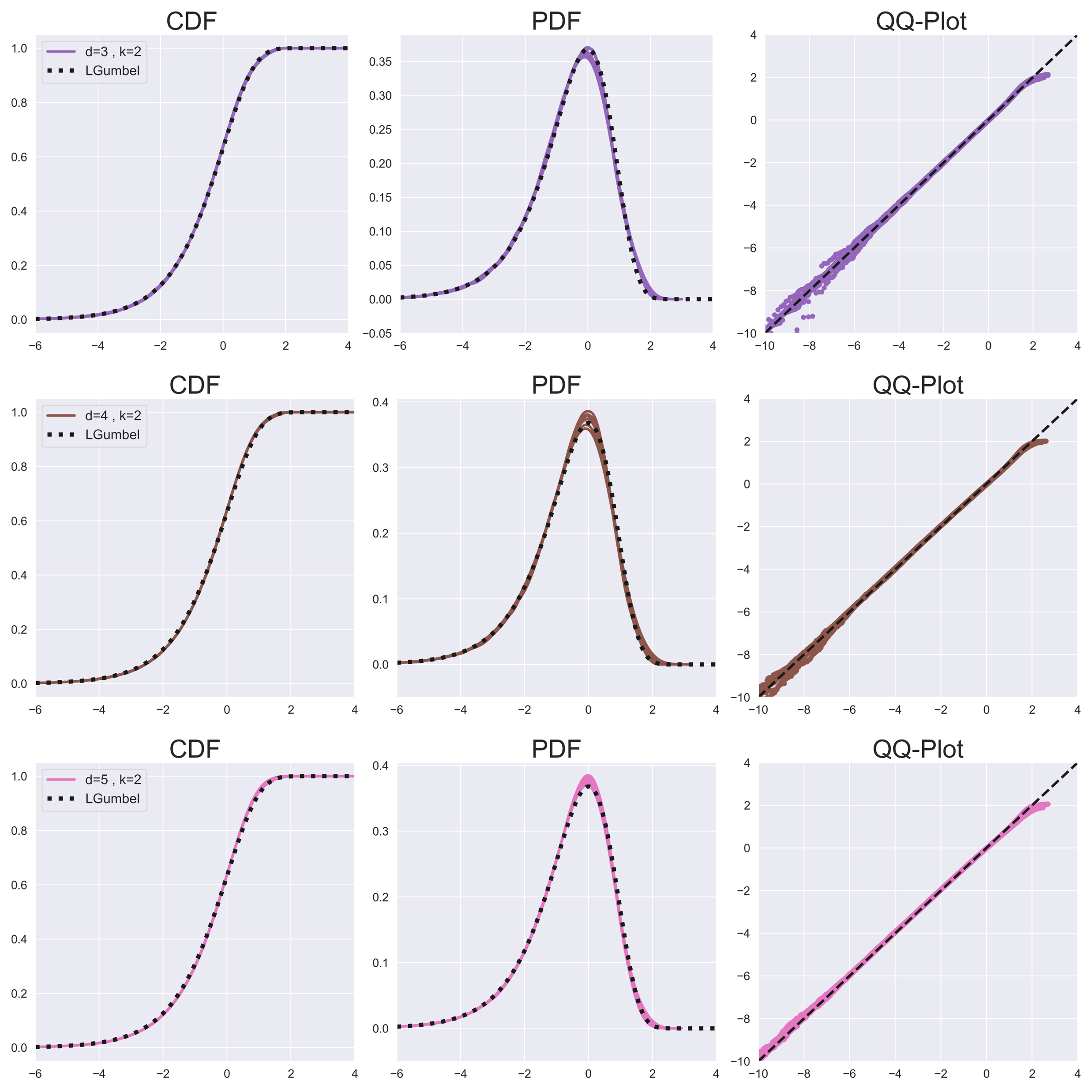}
      \caption{The $\ell$-values distribution for the \v Cech complex -- part 2/3.
    \label{fig:cech_iid_loglog_2}}
 \end{figure}

\begin{figure}[h!]
     \centering
      \includegraphics[width=0.95\textwidth]{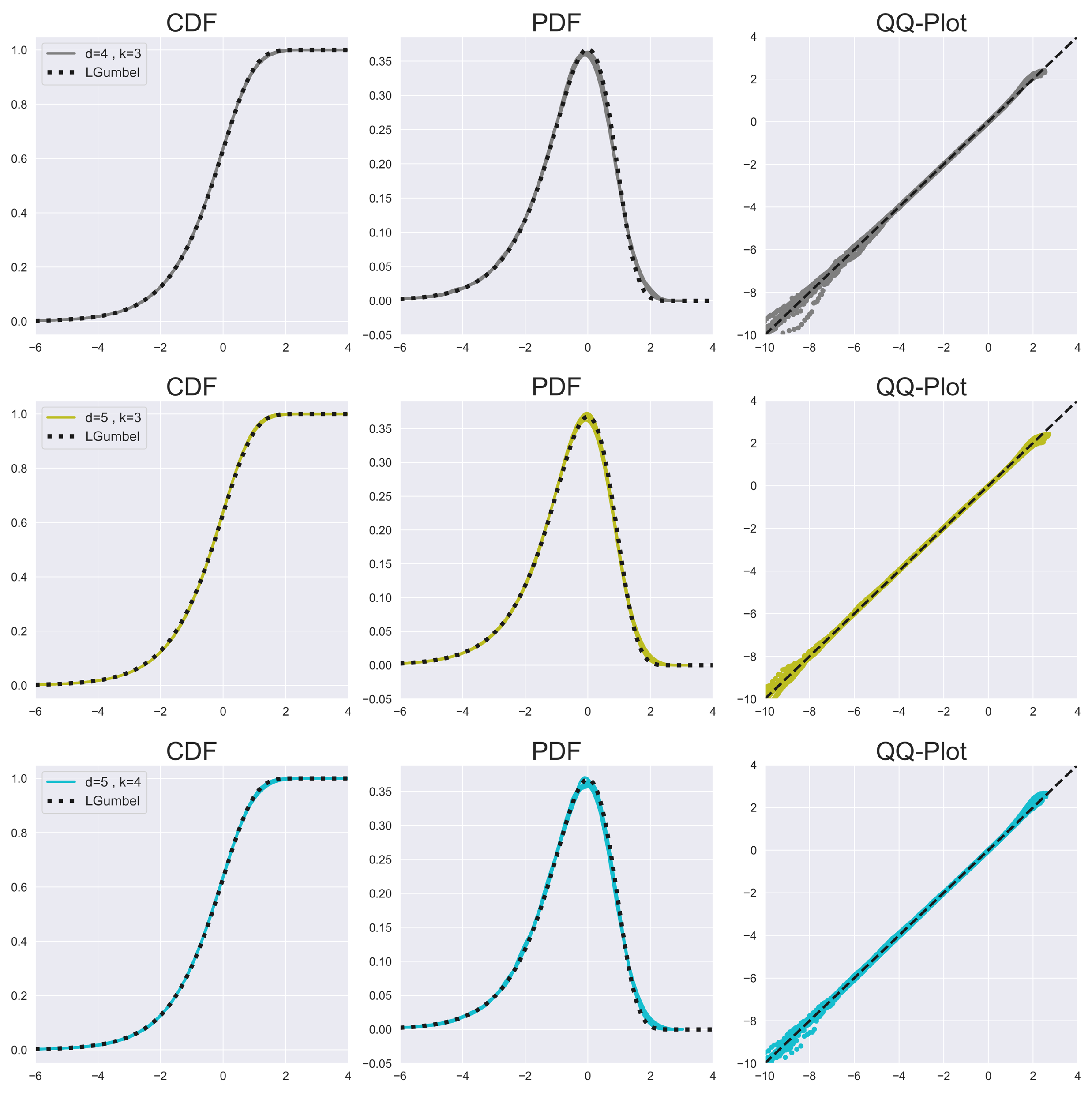}
      \caption{The $\ell$-values distribution for the \v Cech complex -- part 3/3.
    \label{fig:cech_iid_loglog_3}}
 \end{figure}

\begin{figure}[h!]
     \centering
      \includegraphics[width=0.95\textwidth]{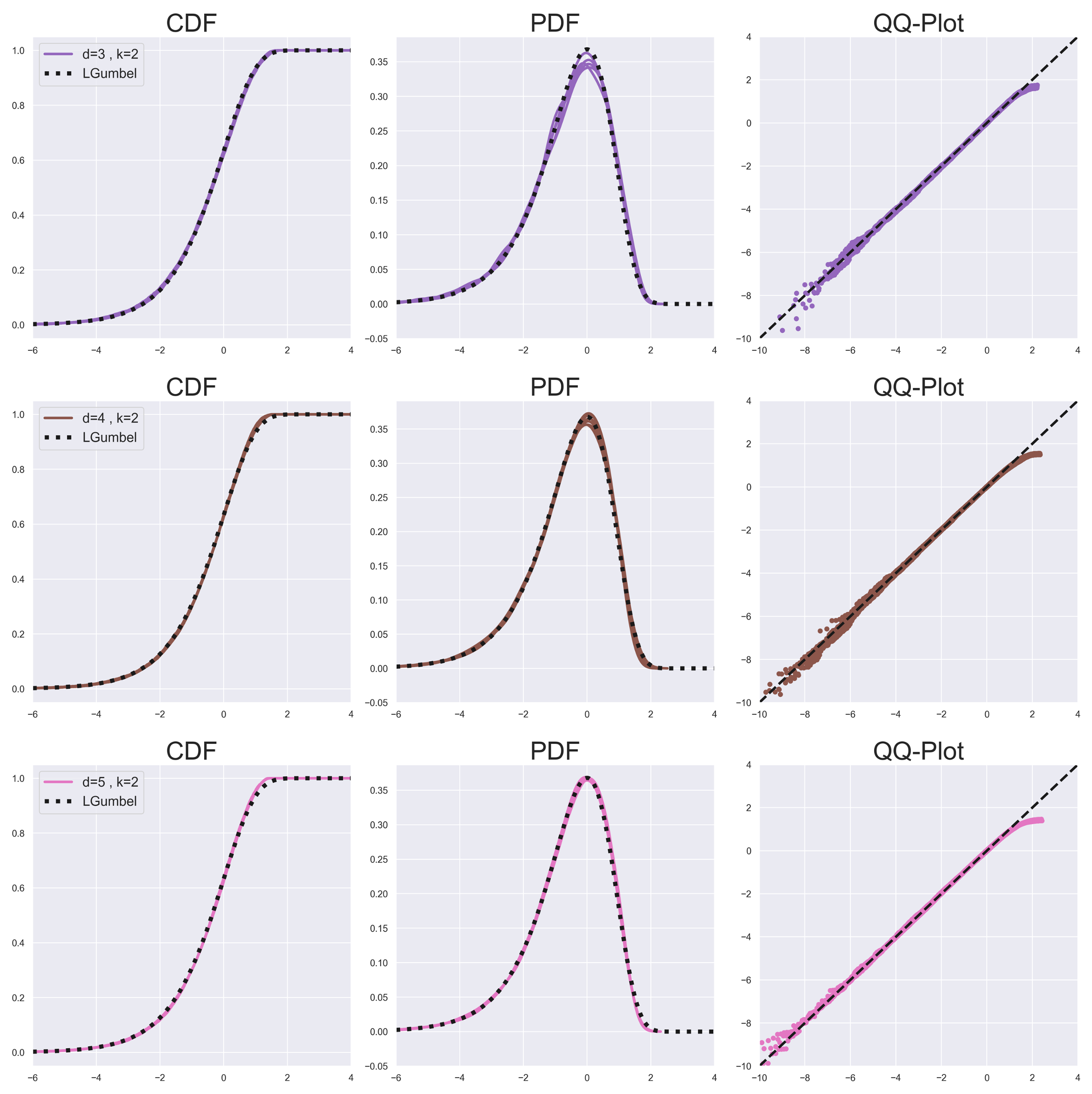}
       \caption{The $\ell$-values distribution for the Rips complex -- part 2/2.
    \label{fig:rips_iid_loglog_2}}
 \end{figure}

We first examine  the experimental results supporting Conjecture \ref{con:iid} more closely. To this end we examined the settings detailed in Tables \ref{tab:iid_uniform} and \ref{tab:iid_nonuniform}. The $\pi$-values distributions for the \v Cech filtration are presented in Figure \ref{fig:cech_iid} and for the Rips filtration in Figure \ref{fig:rips_iid}. We note that not all settings in Tables  \ref{tab:iid_uniform} and \ref{tab:iid_nonuniform} were computed for both the \v Cech and the Rips complex. The reasons for that are purely computational. Some of the settings required more memory or computing power than we have at our disposal. Tables \ref{tab:cech_iid_all}-\ref{tab:rips_iid_all} provide the full list of iid-configurations tested, including the sample size used.
The varying sample sizes are due to both our wish to explore the behavior at different values of $n$, and the computational limitation we have for computing the persistence diagrams.
In Figures \ref{fig:cech_iid_loglog_1}-\ref{fig:rips_iid_loglog_2} we present the distribution of the $\ell$-values for all the iid configurations tested.

\begin{figure}[h!]
     \centering
     \includegraphics[width=\textwidth]{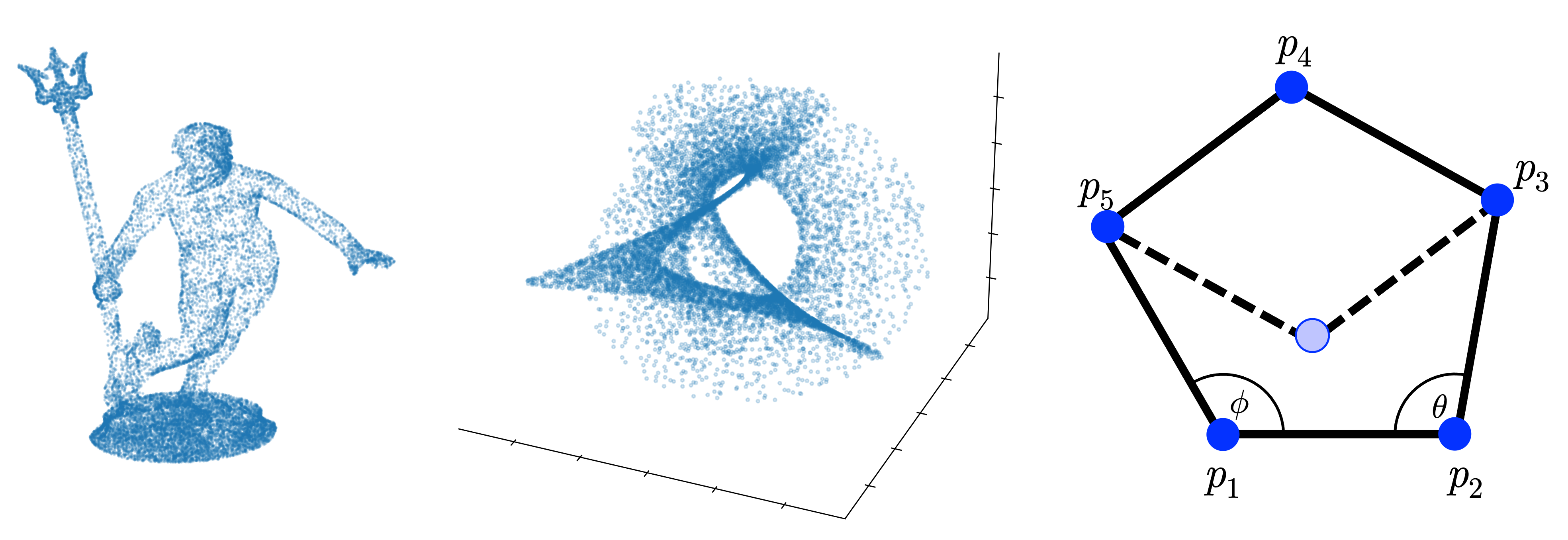}      
     \caption{(left) a uniform sample of the Neptune statue. (middle) a sampling of the Hennenberg surface. (right) A example of a linkage configuration. The light colored vertex is the reflected option mentioned in our description.}
     \label{fig:linkage}
 \end{figure}

Next, we provide mode details about the iid distributions studied.
Let $X = (X^{(1)},\ldots, X^{(d)})$ be a single point in one of the iid samples. We describe briefly how $X$ is generated in each of the models.
\begin{itemize}
    \item {\bf Box:} $X$ is uniformly distributed in $[0,1]^d$.
    \item {\bf Ball:} $X$ is uniformly distributed in a unit $d$-dimensional ball. Sampling was done using the rejection-sampling method.
    \item {\bf Annulus:} $X$ is uniformly distributed in a $d$-dimensional annulus with radii in the range $[1/2,1]$. Sampling was done using the rejection-sampling method.
    \item {\bf Sphere:} $X$ is uniformly distributed on a $d$-dimensional unit sphere, embedded in $\R^{d+1}$. Sampling was done by generating a standard ($d+1$)-dimensional normal variable, and projecting it on the unit sphere.
    \item {\bf Beta(a,b):} The coordinates $X^{(1)},\ldots, X^{(d)}$ are sampled independently from the Beta(a,b) distribution.
    \item {\bf Cauchy:} The coordinates $X^{(1)},\ldots, X^{(d)}$ are sampled independently from the Cauchy distribution.
    \item {\bf Normal:} The coordinates $X^{(1)},\ldots, X^{(d)}$ are sampled independently from the standard normal  distribution.
    \item {\bf Torus:} We generate points on the $2$-dimensional torus embedded in $\R^3$ as follows. We generate two independent variables $\phi$ and $\theta$ uniformly in $[0,2\pi]$. Then we take the coordinates to be:
    \[
        X^{(1)} = (R_1+R_2\cos(\phi))\cos(\theta), \quad
        X^{(2)} = (R_1+R_2\cos(\phi))\sin(\theta),\quad
        X^{(3)} = R_2\sin(\phi).
    \]
    We used $R_1=2$ and $R_1=1$.
    \item {\bf Klein:} We sample an embedding of the Klein bottle into $\mathbb{R}^4$ as follows.
    We generate two independent variables $\phi$ and $\theta$ uniformly in $[0,2\pi]$. The value of $X$ is then computed as
    \begin{align*}
      X^{(1)} &=  (1+\cos(\theta))\cos(\phi)  &&
      X^{(2)} =  (1+\cos(\theta))\sin(\phi)  \\
      X^{(3)} &=  \sin(\theta)\cos\left(\frac{\phi}{2}\right)  &&
      X^{(4)} =  \sin(\theta)\sin\left(\frac{\phi}{2}\right)  
      \end{align*}
    \item {\bf Projective:} We sample an embedding of the real projective plane into $\mathbb{R}^4$. We generate independent variables $U_0$, $V_0$ and $W_0$  from the standard normal distribution.
    Next we take $(U, V, W) = \frac{(U_0,V_0,W_0)}{\sqrt{U_0^2+V_0^2+W_0^2}}\in \S^2$, and define
    \begin{align*}
      X^{(1)} = UV  &&
      X^{(2)} = UW &&  
      X^{(3)} = V^2 - W^2&&
      X^{(4)} = 2VW.
      \end{align*}
    \item {\bf Linkage:} This model samples the configuration space of  unit pentagonal linkages, i.e. a pentagon where adjacent edges are of unit length. To sample this space, we first fix two vertices at $p_1=(0,0)$ and $p_2=(1,0)$ respectively. Next, we generate two independent variables $\phi$ and $\theta$ uniformly in $[0,2\pi]$, and set \begin{align*}
       p_5=(\cos(\varphi),\sin(\varphi)), && p_3=(1+\cos(\theta),\sin(\theta)). 
    \end{align*}
    If $\|p_3-p_5\|>2$, the sample is rejected as there is no linkage with the chosen angles. Otherwise, there are two possible choices of the last point $p_4$. Let $(q^{(1)},q^{(2)})$ denote the midpoint of $\overline{p_3 p_5}$. Then 
    $$
    p_4 = q + S \frac{\sqrt{1-\|q-p_5\|^2}}{\|q-p_5\|}(p_5^{(2)} -q^{(2)},q^{(1)}-p_5^{(2)}) 
    $$
    where $S$ is independent of $\phi,\theta$, and
    $\mathbb{P}(S=1) = \frac{1}{2}$ and $\mathbb{P}(S=-1) = \frac{1}{2} $. See Figure \ref{fig:linkage} for an example.
    \item {\bf Neptune:} We construct a sample from the surface of the statue of Neptune.  From \cite{aimshape}, we retrieved a triangulation of the surface, consisting of 4,007,872 triangles. To generate a sample, we first compute the area of each triangle and then choose a triangle at random with probability inversely proportional to the area of the triangle. We then pick a point uniformly in the chosen triangle. 
    \item {\bf Hennenberg:} We construct a sample of the Henneberg surface in $\mathbb{R}^3$.
    We start by generating two independent variables $\phi$ and $\theta$ uniformly in $[0,2\pi]$.
    We then construct the sample by
    \begin{align*}
        X^{(1)} &= 2\cos(\theta)\sinh(\phi) - \frac{2}{3}\cos(3\theta)\sinh(3\phi),\\
        X^{(2)} &= 2\sin(\theta)\sinh(\phi) + \frac{2}{3}\sin(3\theta)\sinh(3\phi),\\
        X^{(3)} &= 2\cos(2\theta)\cosh(2\phi).
    \end{align*}
    \item {\bf Stratified Spaces:} To construct $X$, we consider two spaces $M_1 \subset M_2$ such that the dimension of $M_1$ is less than $M_2$. Then one of the two spaces is chosen with some probability $p$ (which is a parameter of the model), and the chosen space is sampled uniformly. Although several models were tried, in the examples we show a plane $[-1,1]^2$ embedded in the middle of a cube $[-1,1]^3$. In other words, if the point $(x,y)$ is chosen from the plane, the coordinates would be  $(x,y,0)$. 
\end{itemize}
 
\begin{table}[h!]
    \centering
    \begin{tabular}{|c|c|c|c|}
    \hline
        distribution ($\P$)&
        dimension ($d$)&
        homology degree ($k$)&
        sample size ($n$) \\
        \hline
        \hline
        		Box  & 2 & 1& 1,000,000\\
		Box  & 3 & 1-2& 1,000,000\\
		Box  & 4 & 1-3 & 50,000\\
		Box  & 5 & 1-4 &100,000\\
		Ball  & 2 & 1 & 50,000\\
		Ball  & 3 & 1-2 & 50,000\\
		Ball & 4 & 1-3 & 50,000\\
		Ball  & 5 & 1-4 & 50,000\\
		Annulus  & 2 & 1 & 50,000 \\	
		Annulus & 3 & 1-2 & 10,000 \\	
		Annulus & 4 & 1-3 & 10,000 \\	
		Annulus & 5 & 1-4 & 10,000 \\	
		Sphere  & 2 & 1 & 100,000 \\	
		Sphere & 3 & 1-2 & 50,000 \\	
		Sphere & 4 & 1-3 & 50,000 \\	
		Sphere & 5 & 1-4 & 10,000 \\	
		Beta(3,1) & 2 & 1 & 50,000\\
		Beta(3,1) & 3 & 1-2 & 50,000\\
		Beta(3,1) & 4 & 1-3 & 50,000\\
		Beta(3,1) & 5 & 1-4 & 10,000\\
		Beta(5,2) & 2 & 1 & 100,000\\
		Beta(5,2) & 3 & 1-2 & 50,000\\
		Beta(5,2) & 4 & 1-3 & 50,000\\
		Beta(5,2) & 5 & 1-4 & 10,000\\
		Cauchy & 2 & 1 & 100,000\\
		Cauchy & 3 & 1-2 & 100,000\\
		Cauchy & 4 & 1-3 & 100,000\\
		Cauchy & 5 & 1-4 & 10,000\\
		Normal & 2 & 1 & 1,000,000\\
		Normal & 3 & 1-2 & 1,000,000\\
		Normal & 4 & 1-3 & 50,000\\
		Normal & 5 & 1-4 & 50,000\\
		Torus & 2 & 1 & 1,000,000\\
		Klein & 2 & 1 & 1,000,000\\
		Projective & 2 & 1 & 50,000\\
                         \hline
    \end{tabular}
    \caption{Complete experimental list for the  iid settings used with the \v Cech complex.
    \label{tab:cech_iid_all}}
\end{table}

\begin{table}[h!]
    \centering
    \begin{tabular}{|c|c|c|c|}
    \hline
        distribution ($\P$)&
        dimension ($d$)&
        homology degree ($k$)&
        sample size ($n$) \\
        \hline
        \hline
        		Box  & 2 & 1& 100,000\\
		Box  & 3 & 1-2& 50,000\\
		Box  & 4 & 1-2 & 50,000\\
		Box  & 5 & 1-2 & 50,000\\
		Ball  & 2 & 1 & 50,000\\
		Ball  & 3 & 1-2 & 50,000\\
		Ball & 4 & 1-2 & 50,000\\
		Ball  & 5 & 1-2 & 50,000\\
		Annulus  & 2 & 1 & 10,000 \\	
		Annulus & 3 & 1-2 & 50,000 \\	
		Annulus & 4 & 1-2 & 10,000 \\	
		Annulus & 5 & 1-2 & 10,000 \\	
		Sphere  & 2 & 1 & 50,000 \\	
		Sphere & 3 & 1-2 & 50,000 \\	
		Sphere & 4 & 1-3 & 100,000 \\	
		Sphere & 5 & 1-4 & 50,000 \\	
		Beta(3,1) & 2 & 1 & 50,000\\
		Beta(3,1) & 3 & 1-2 & 50,000\\
		Beta(3,1) & 4 & 1-2 & 50,000\\
		Beta(3,1) & 5 & 1-2 & 10,000\\
		Beta(5,2) & 2 & 1 & 50,000\\
		Beta(5,2) & 3 & 1-2 & 50,000\\
		Beta(5,2) & 4 & 1-2 & 50,000\\
		Beta(5,2) & 5 & 1 & 10,000\\
		Cauchy & 2 & 1 & 10,000\\
		Cauchy & 3 & 1-2 & 50,000\\
		Cauchy & 4 & 1-2 & 50,000\\
		Cauchy & 5 & 1-2 & 10,000\\
		Normal & 2 & 1 & 50,000\\
		Normal & 3 & 1-2 & 50,000\\
		Normal & 4 & 1-2 & 50,000\\
		Normal & 5 & 1-2 & 50,000\\
		Torus & 2 & 1 & 1,000,000\\
		Klein & 2 & 1 & 10,000\\
		Projective & 2 & 1 & 50,000\\
		Linkage & 2 & 1 & 50,000\\
		Neptune & 2 & 1 & 50,000\\
                         \hline
    \end{tabular}
    \caption{Complete experimental list for the iid settings used with the Rips complex.
    \label{tab:rips_iid_all}}
\end{table}
 
 \newpage 
\FloatBarrier
 \subsection{Sampling from non-iid distributions}\label{app:non-iid}

 As mentioned in the paper, we tested two examples of non-iid point-clouds.
 The first example is sampling the path of a $d$-dimensional Brownian motion $W_t$. To sample $W_t$ at times $t=1,\ldots,n$ we use the fact that $W_t$ has stationary independent increments. We start by taking $Z_1,\ldots,Z_n$ to be iid $d$-dimensional standard normal variables, and then we define $X_i = W_{t=i}= \sum_{i=1}^n Z_i$. 
 
\begin{figure}[h!]
     \centering
      \includegraphics[width=0.95\textwidth]{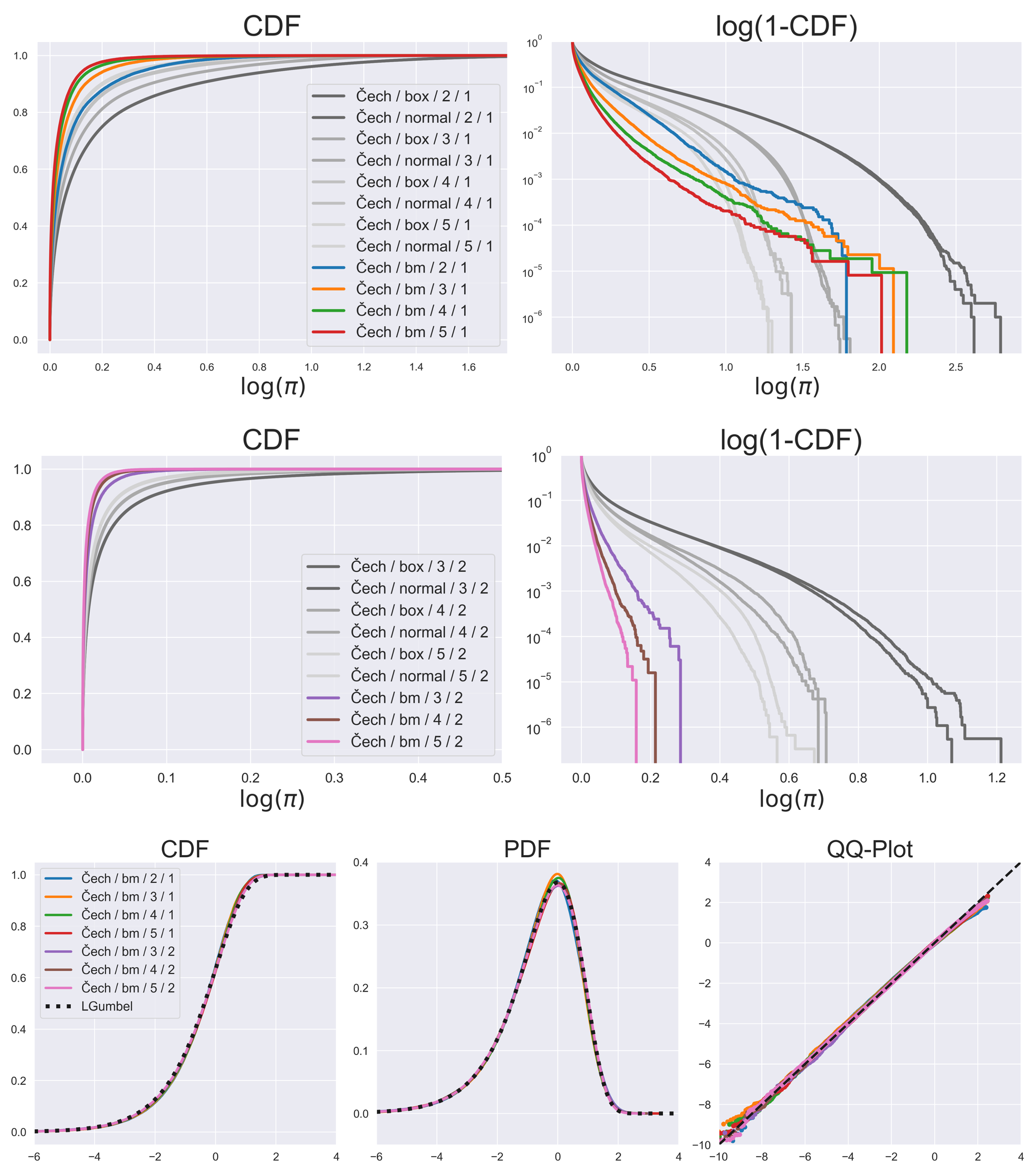}
      \caption{Sampling from a Brownian motion, the \v Cech complex.
    \label{fig:cech_bm}}
 \end{figure}

\begin{figure}[h!]
     \centering
      \includegraphics[width=0.95\textwidth]{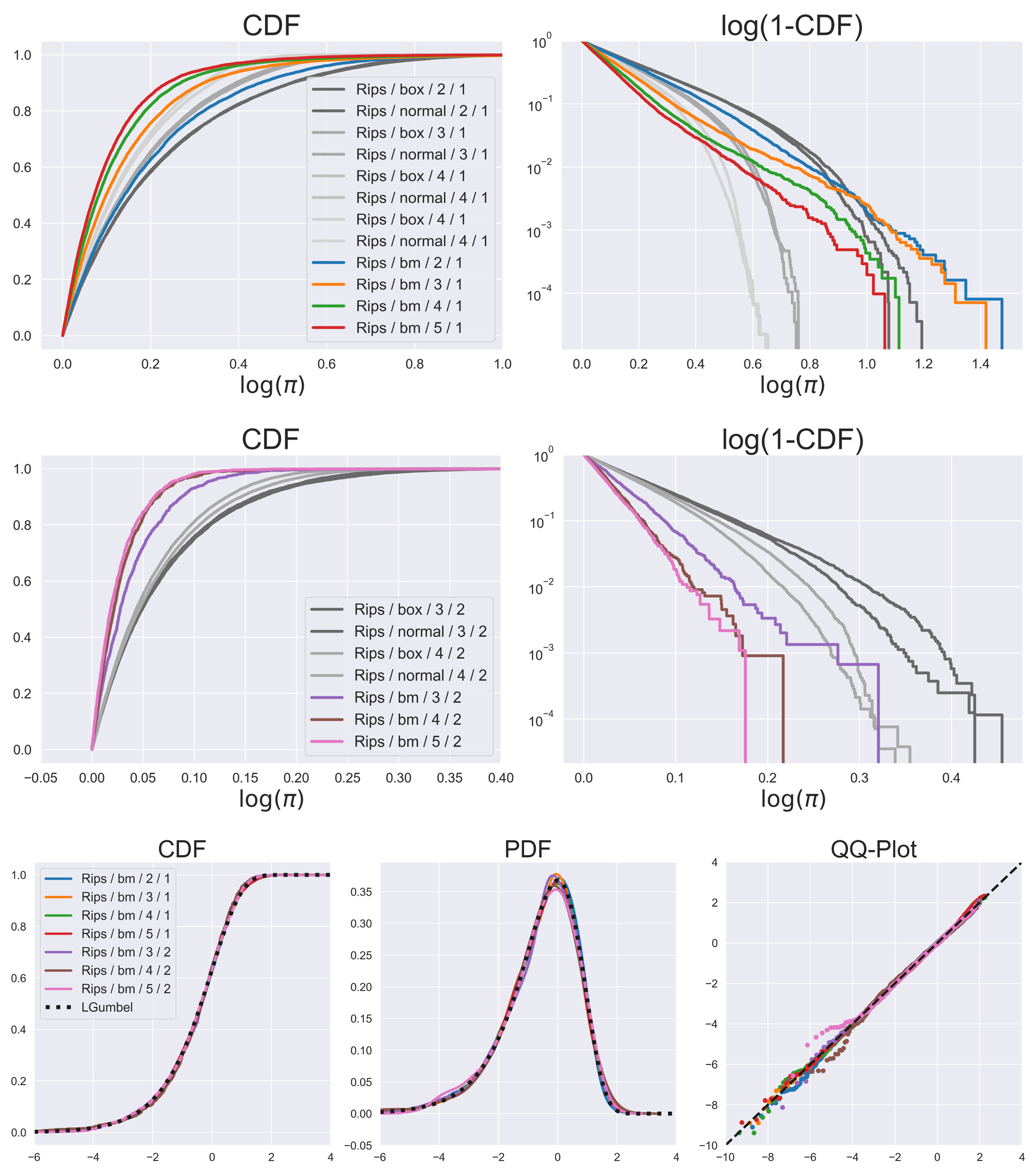}
      \caption{Sampling from a Brownian motion, the Rips complex.
    \label{fig:rips_bm}}
 \end{figure}

 In the top two rows of Figures \ref{fig:cech_bm}-\ref{fig:rips_bm} we present the distribution of $\pi$-values for Brownian point-clouds (in color) compared to a selected baseline from the iid distributions (in gray). We included these figures here to show that the Brownian  samples behave very differently from the iid samples. Nevertheless, the bottom two row presents the distribution of the $\ell$-values, which seems to follow the same LGumbel distribution as the iid case. All samples here are of size $n=50,000$.

\begin{figure}[h!]
     \centering
      \includegraphics[width=0.95\textwidth]{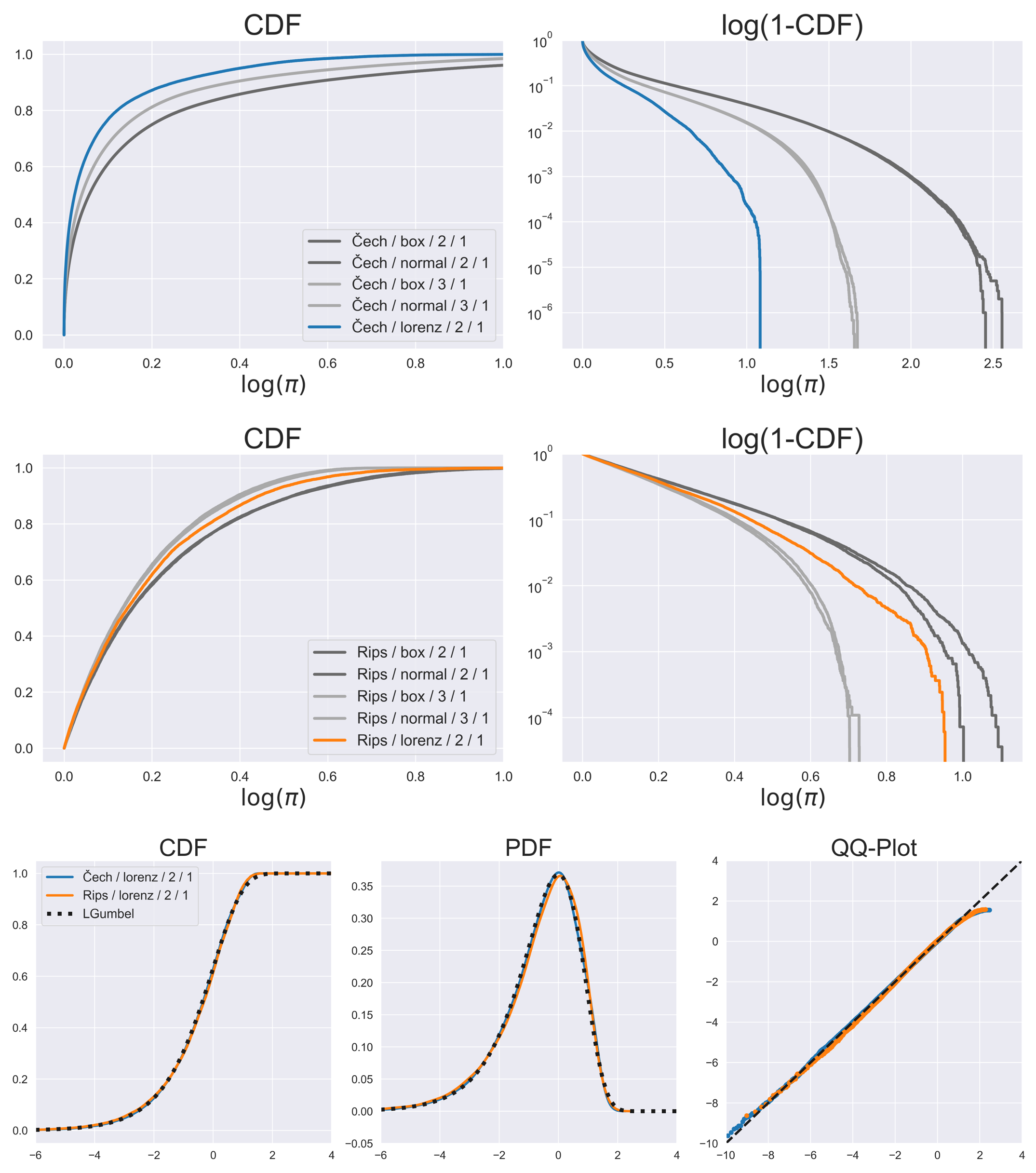}
      \caption{Sampling the Lorenz system -- both \v Cech and Rips.
    \label{fig:lorenz}}
 \end{figure}

 The second example we examined is a discrete-time sample of the Lorenz dynamical system, which is generated as follows. We start by picking a random initial point, uniformly in $[0,1]^3$. We then generate the sample using the differential equations,
    \begin{align*}
        \frac{dX^{(1)}}{dt} &= \sigma \left(X^{(2)} - X^{(1)}\right),\\
        \frac{dX^{(2)}}{dt} &= X^{(1)} \left(\rho - X^{(3)}\right)- X^{(2)},\\
        \frac{dX^{(3)}}{dt} &= X^{(1)}X^{(2)} - \beta X^{(3)}.
    \end{align*}
    We use $\sigma=45$, $\rho = 54$, and $\beta=10$, and a numerical approximation with $dt=0.1$ to generate a trajectory for the number of samples required. Note that each instance is a single trajectory. The results are presented in Figure \ref{fig:lorenz}. The sample size for the \v Cech is $n=100,000$, and for the Rips is $n=50,000$.
    
 \newpage\FloatBarrier

 \subsection{Testing on real data} \label{app:real-data}
 
 As described in the paper we tested the strong-universality conjectures against two examples of real data - patches of natural images, and sliding windows of voice recordings. In this section we provide more details about these examples.

\noindent {\bf Image patches.} 
 The images were taken from van Hateren and van der Schaaf \cite{van1998independent} database. This database contains a collection of  about 4,000 gray-scale images.
 We follow the procedure described in \cite{carlsson_local_2008}.
 We randomly select patches of size $3\times 3$ from the entire dataset. This gives us a point-cloud in $\R^9$. Let $\bx = (x_1,\ldots,x_9)\in \R^9$ represent the log-values of a single patch, and follow the following steps:
 \begin{enumerate}
     \item Compute the average pixel value $\bar \bx$, and subtract it from all 9 pixels in the patch, i.e. $\by = \bx-\bar\bx$
     \item Compute the ``D-norm'' (a measure for contrast), $\|\by\|_D = \sqrt{\by^TD\by}$. 
     \item Use the D-norm to normalize the pixel values, $\bz = \by / \norm{\by}_D$.
     \item Use the Discrete Cosine Transform (DCT) basis to change the coordinate system, $\bv = \Lambda A^T \bz$.
 \end{enumerate}
 The values of $D,A,\Lambda$, as well as more details about this proecedure can be found in \cite{lee_nonlinear_2003,carlsson_local_2008}.
 The process above results in a point-cloud lying on the unit 7-dimensional sphere in $\R^8$. 
 
 For the topological analysis in \cite{carlsson_local_2008}, the point-cloud is further filtered, in order to focus on the ``essential'' information captured by the patches. This is done in two steps:
 \begin{enumerate}
     \item Keep only  ``high-contrast'' patches -- whose D-norm is in the top 20\%.
     \item Compute the distance of each of the remaining patches to their k-nearest neighbor (with $k=15$), and keep only the patches in the bottom 15\%.
 \end{enumerate}
 We repeat the exact procedure performed in \cite{carlsson_local_2008} for two main reasons. Firstly, in Section \ref{sec:ht_fin} we wish to use our hypothesis testing framework to assign p-values to the cycles found using this procedure. Secondly, this procedure makes the sample distribution more intricate, and adds dependency between the sample points. We wanted to challenge our conjectures with data as complex as possible.
 
\begin{figure}[h!]
     \centering
      \includegraphics[width=0.95\textwidth]{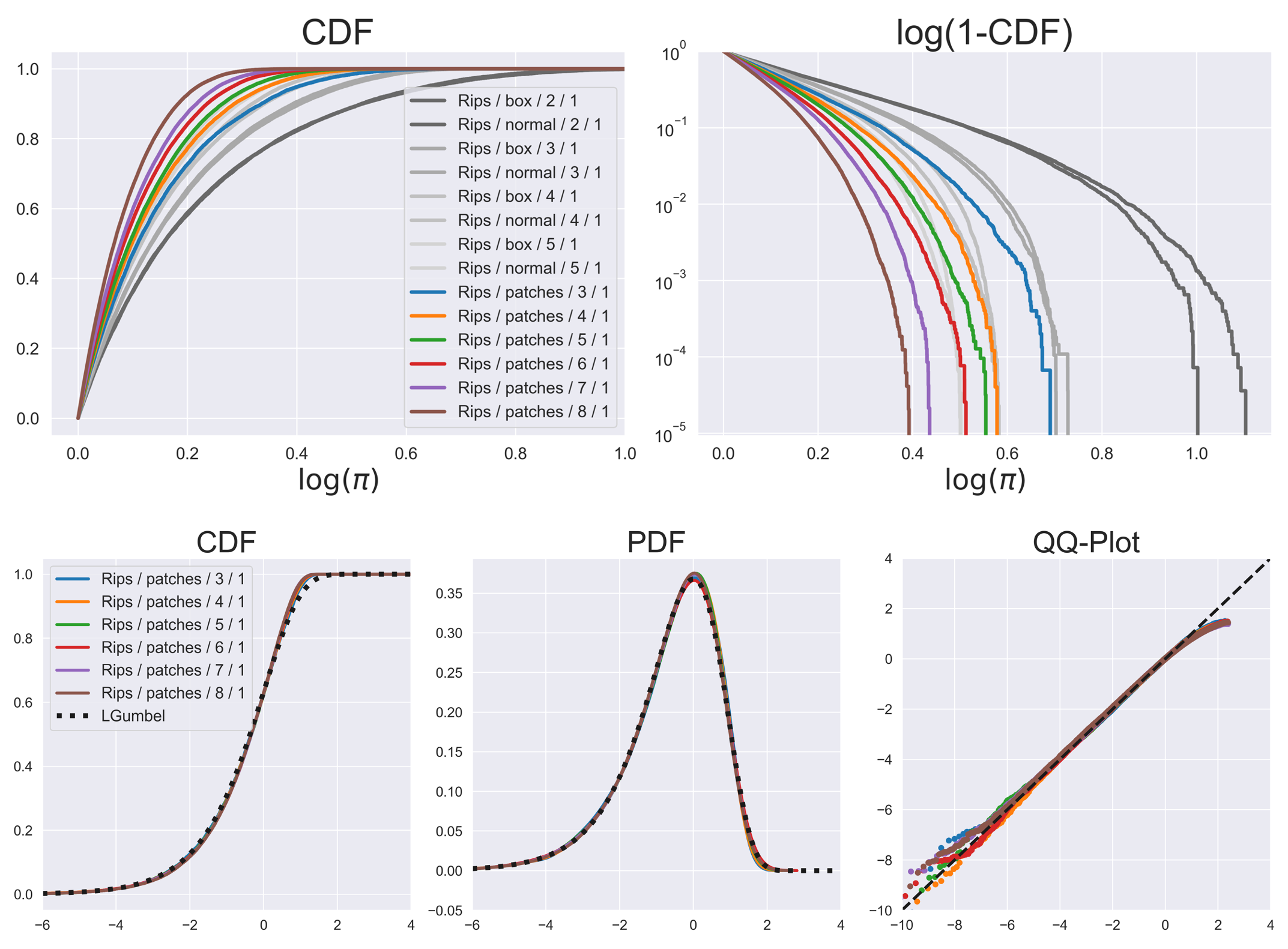}
      \caption{Sampling patches of natural image, Rips complex.
    \label{fig:patches}}
 \end{figure}

 The results for this point-cloud are presented in Figure \ref{fig:patches}. In addition to taking the original 8-dimensional point-cloud, we also examined its lower dimensional projections for dimensions $d=3,\ldots,7$. The sample size used was $n=50,000$ for all dimensions.

\begin{figure}[h!]
     \centering
      \includegraphics[width=0.95\textwidth]{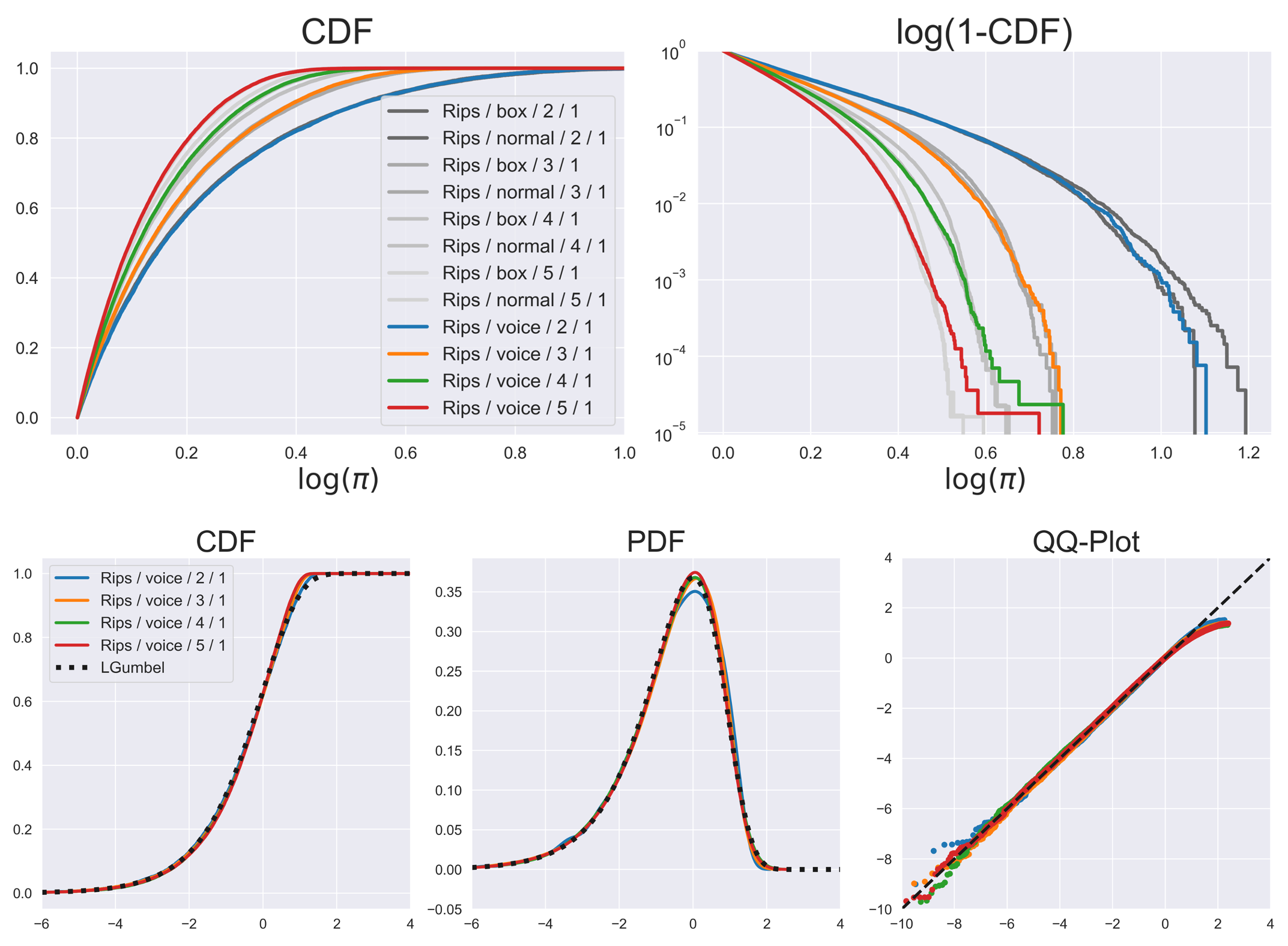}
      \caption{Sampling sliding windows from a voice recording, Rips complex.
    \label{fig:voice}}
 \end{figure}

\noindent {\bf Sound recordings.} 
We took an arbitrary audio recording (the voice of one of the authors), and applied the time-delay embedding transformation (cf.\cite{takens_detecting_1981}) to convert the temporal signal into a $d$-dimensional point-cloud.
The voice recording is a 20 seconds excerpt, sampled at 16KHz, 48Kbps. We denote the corresponding discrete time signal as $V_t$. To convert the signal into a point-cloud we used the following:
\[
    X_i = (V_{i\Delta}, V_{i\Delta+\tau},\ldots, V_{i\Delta+(d-1)\tau}) \in \R^d,\quad i=1,\ldots,n.
\]
The values we took for the example here were $\Delta = 3$, and $\tau = 7$. This, in particular, generates overlap between the windows, which guarantees strong dependency both between the points, and between the coordinates of each point. The sample size taken was $n=50,000$.
The results are presented in Figure \ref{fig:voice}.

\section{Dependency testing}\label{app:indep}

In Section \ref{sec:strong}, we provided evidence for the pairwise independence of the different persistence $\ell$-values. We wish to provide more details here. 

\begin{figure}[h!]
    \centering
    \includegraphics[width=0.7\textwidth]{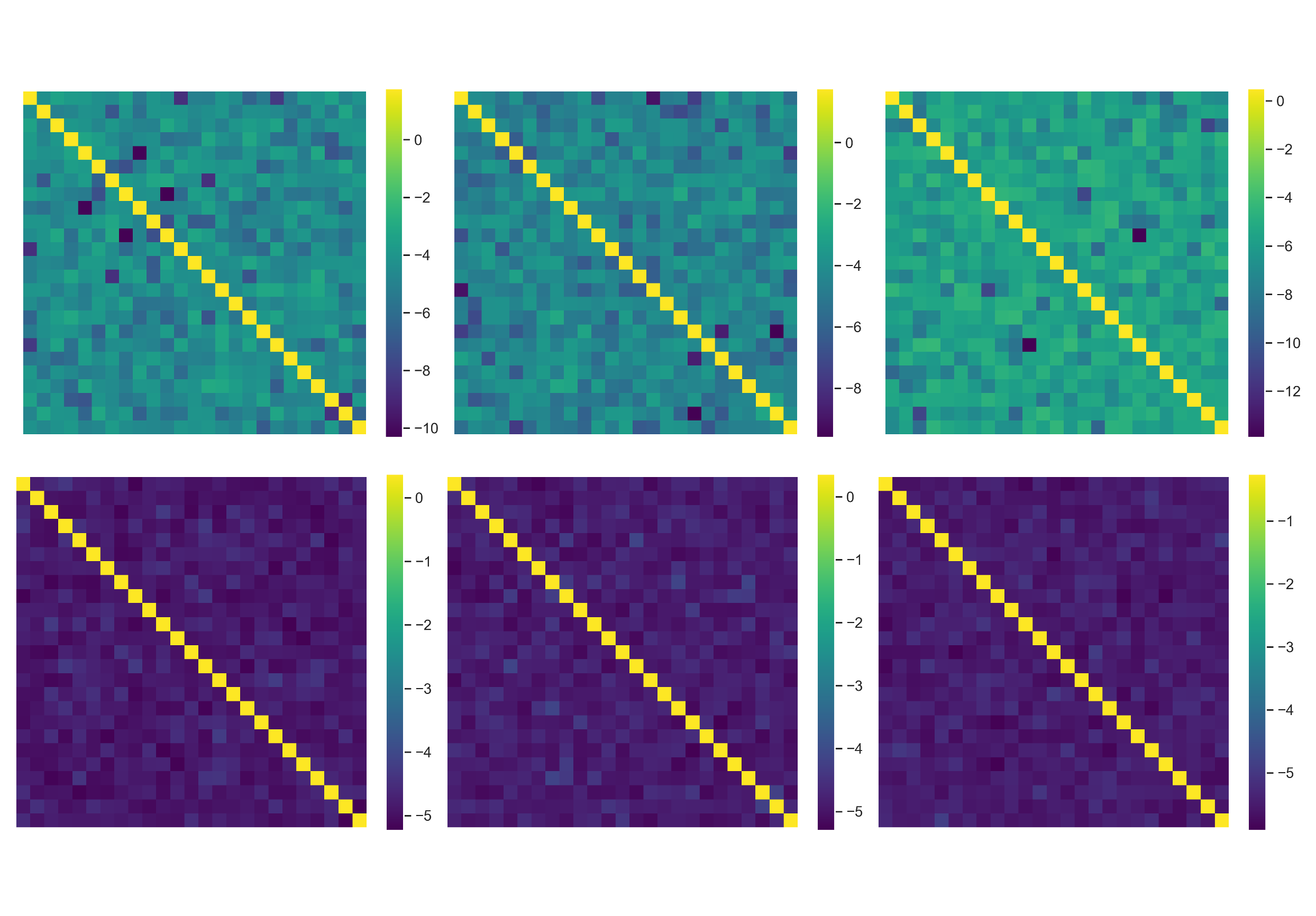}
    \caption{The covariance matrix (top) and distance covariance matrix (bottom), for $\ell_1,\ldots,\ell_{25}$, estimated from $N=100,000$ repeated trials. (left) The $\ell$-values are generated from the uniform distribution on the $2$-dimensional box. (middle) The $\ell$-values are generated by sampling the $2$-dimensional Brownian motion. (right) The $\ell$-values here are simply taken to be 25 iid random LGumbel variables. }
    \label{fig:indep_all}
\end{figure}

As stated in the paper, since points in a persistence diagram have no natural ordering, we used a random ordering on the persistence values. 
Specifically, for each $1\le i \le N$, we generated a new point cloud, computed a persistence diagram, randomly sampled $25$ points from the diagram $p_1,\ldots,p_{25}$, and placed their corresponding $\ell$-values in a vector of the form $L_i = (\ell(p_1),\ldots,\ell(p_{25})) \in \R^{25}$.
We took $N=100,000$ different realizations, from which we estimated both the covariance and the distance-covariance. As these experiments are quite costly (computing $100,000$ persistence diagrams for an individual setting), we only examined two cases here -- the 2-dimensional box, and the 2-dimensional Brownian sample, in order to compare an iid setting with a non-iid one.
The results are presented in Figure \ref{fig:indep_all}. In addition to the results for the box (right), and the Brownian motion (middle), we also included the covariance and distance-covariance matrices, computed from 25 independent LGumbel variables. In addition to the matrices look rather similar, we also include the mean and maximum absolute values form these matrices in Table \ref{tab:indep}.

\begin{table}[]
    \centering
    \begin{tabular}{|c|c|c|c|c|}
    \hline
        & mean-corr & mean d-corr & max-corr & max-dcorr
        \\
        \hline
        \hline
            2d-box & 0.00233 & 0.00535 &  0.00887 & 0.01002 \\
            
            2d-brownian & 0.00248 & 0.00540 & 0.01139 & 0.01128\\
            
            iid-LGumbel & 0.00242 & 0.00525 & 0.00843 &0.00999
            \\
                         \hline
    \end{tabular}
    \caption{Comparing  correlations between $\ell$-values in a box, Brownian motion, and iid LGumbel variables. We compute the mean and maximum absolute values for both the correlation matrix and the distance-correlation of the 25 $\ell$-values extracted from the persistence diagrams.
      \label{tab:indep}}
\end{table}

 \section{Exploring the value of $B$}
The value of the parameter $B$ in \eqref{eqn:loglog} is the only part of our conjectures that depends on  model information. While in our experiments, it seems that we are able to estimate it from the diagram itself quite accurately, we want to further explore its behavior, and in particular how it depends on the model parameters.
As already seen in Figure~\ref{fig:means}, the value of $B$  captures information such as the complex type ($\cT$), dimension ($d$), and homological dimension ($k$). In this section we investigate the relationship further.

\begin{figure}[h]
\begin{center}
		\includegraphics[width=0.5\columnwidth]{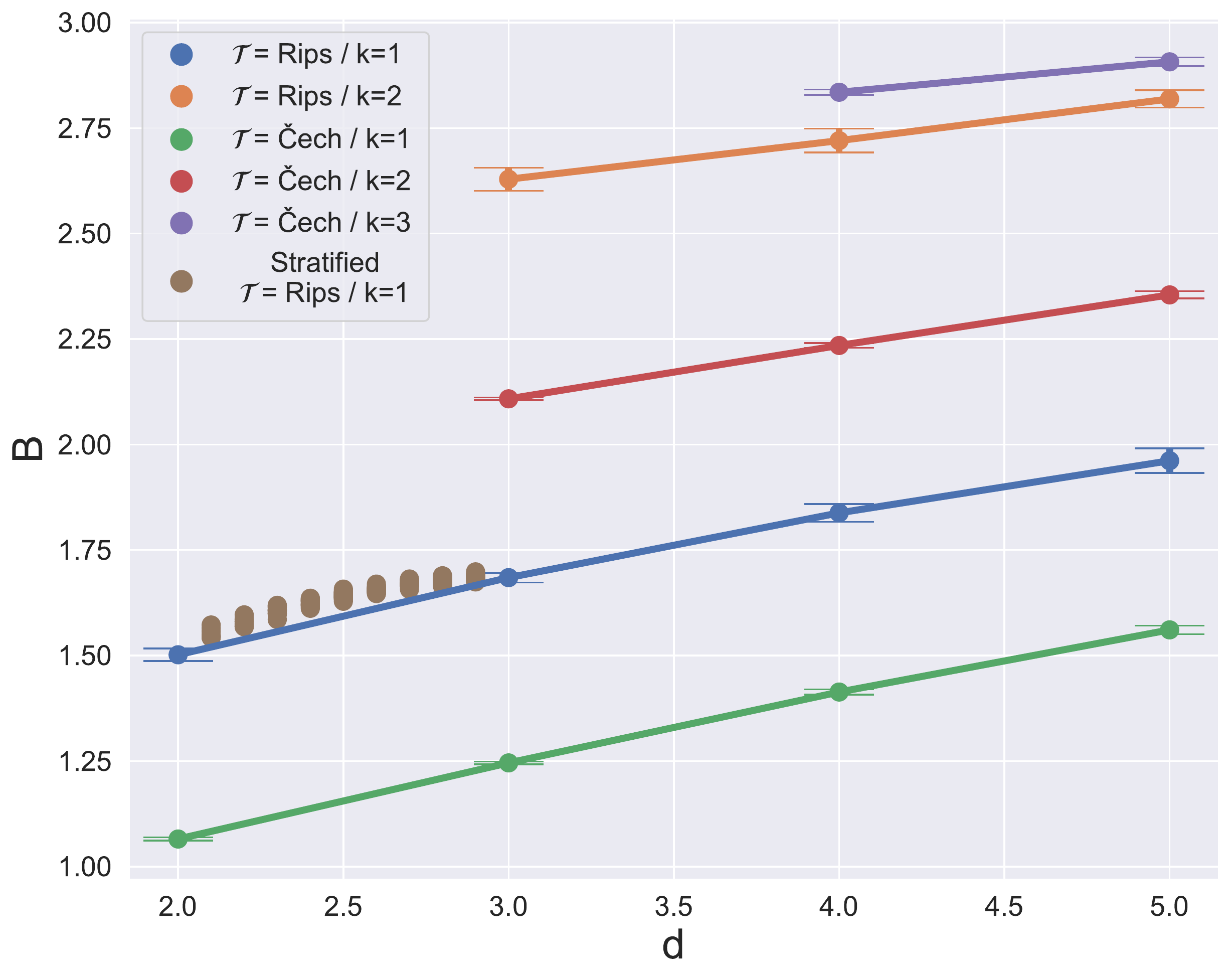}
 	\caption{The estimated value for $B$, for iid samples with  various choices of $(d,\cT,k)$. We present the mean $\pm$standard deviation of our estimates across various models. 
  	Here, we added stratified space samples, which are mixtures of different dimensional spaces (a plane embedded in a cube). 
\label{fig:means_full}}
 	\end{center}
\end{figure}

Figure~\ref{fig:means_full}, is similar Figure~\ref{fig:means}, where we explore the estimated values of $B$, showing the means (with error bars) versus the dimension, for the different complex types and homological dimensions. Here we also include the experiments on  stratified samples, which are iid but of mixed dimensions (between 2 and 3). For the dimension value (the x-axis), we use the weighted average corresponding to the percentage of points sampled from each  stratum. We can see that the values are interpolating between dimension 2 and 3, but this is not a linear relationship. 

\begin{figure}[h!]
     \centering
     \begin{subfigure}{\textwidth}
       \includegraphics[width=0.95\textwidth]{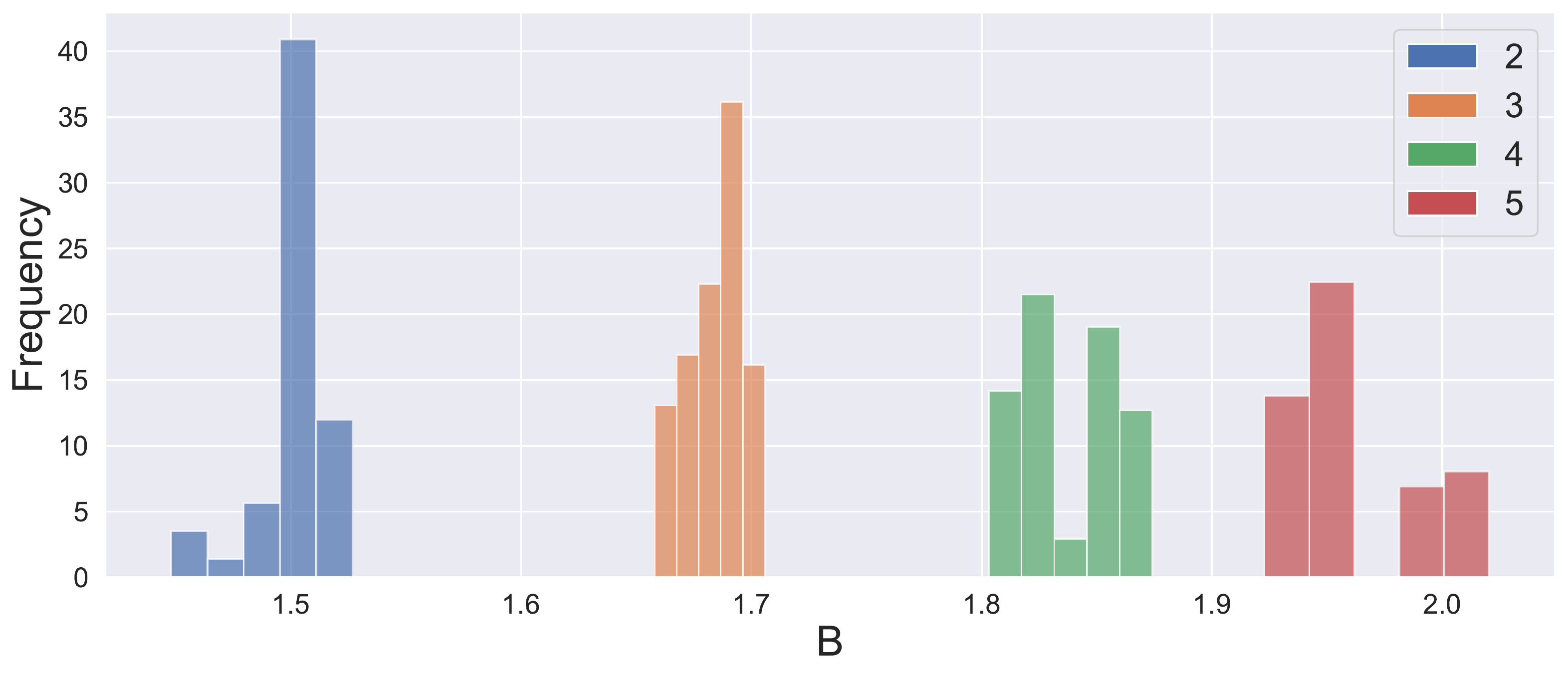}
     \end{subfigure}
     \begin{subfigure}{\textwidth}
        \includegraphics[width=0.95\textwidth]{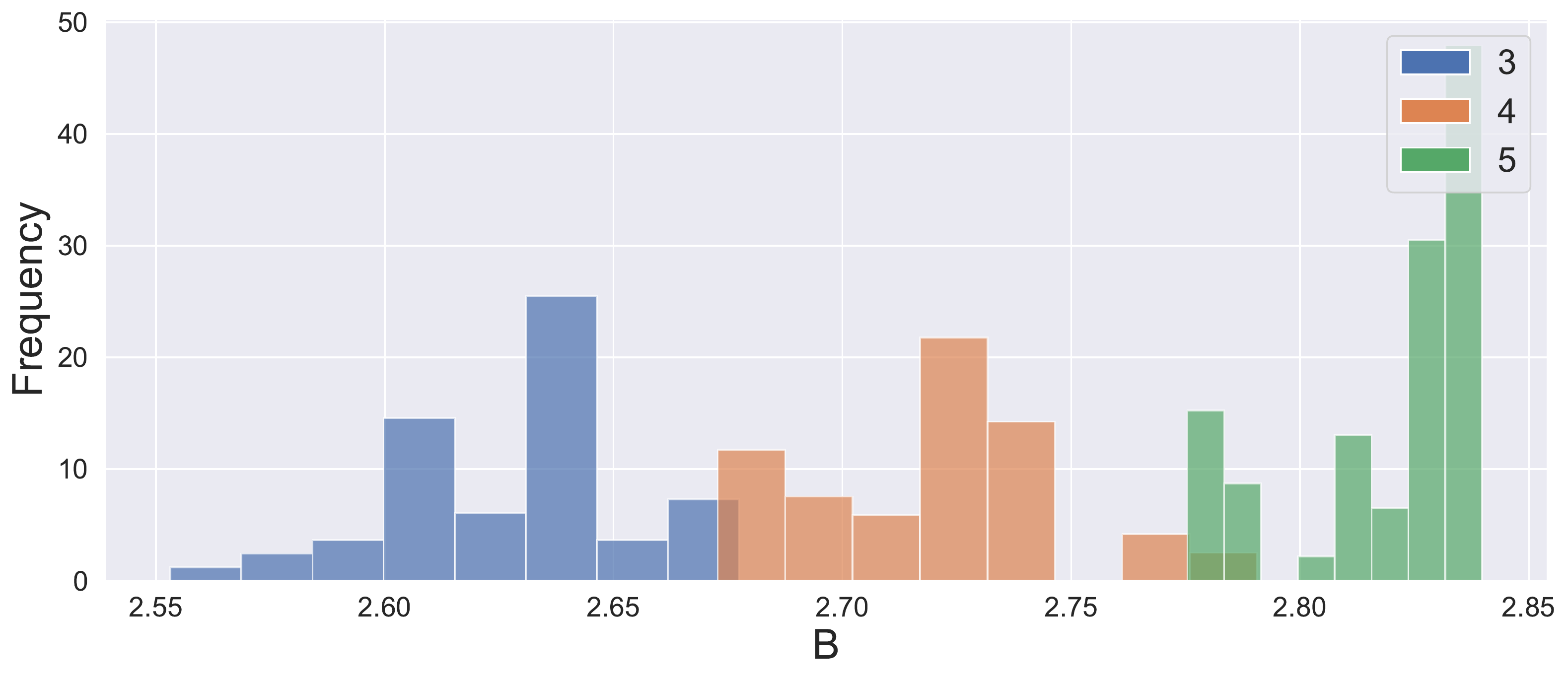}
     \end{subfigure}
       \caption{
       The distributions of the estimated $B$ values, for the Rips complex, across different models for $k=1$ (top) and $k=2$ (bottom).
    \label{fig:Rips_H1_B}}
 \end{figure}

\begin{figure}[h!]
     \centering
     \begin{subfigure}{\textwidth}
         \includegraphics[width=0.95\textwidth]{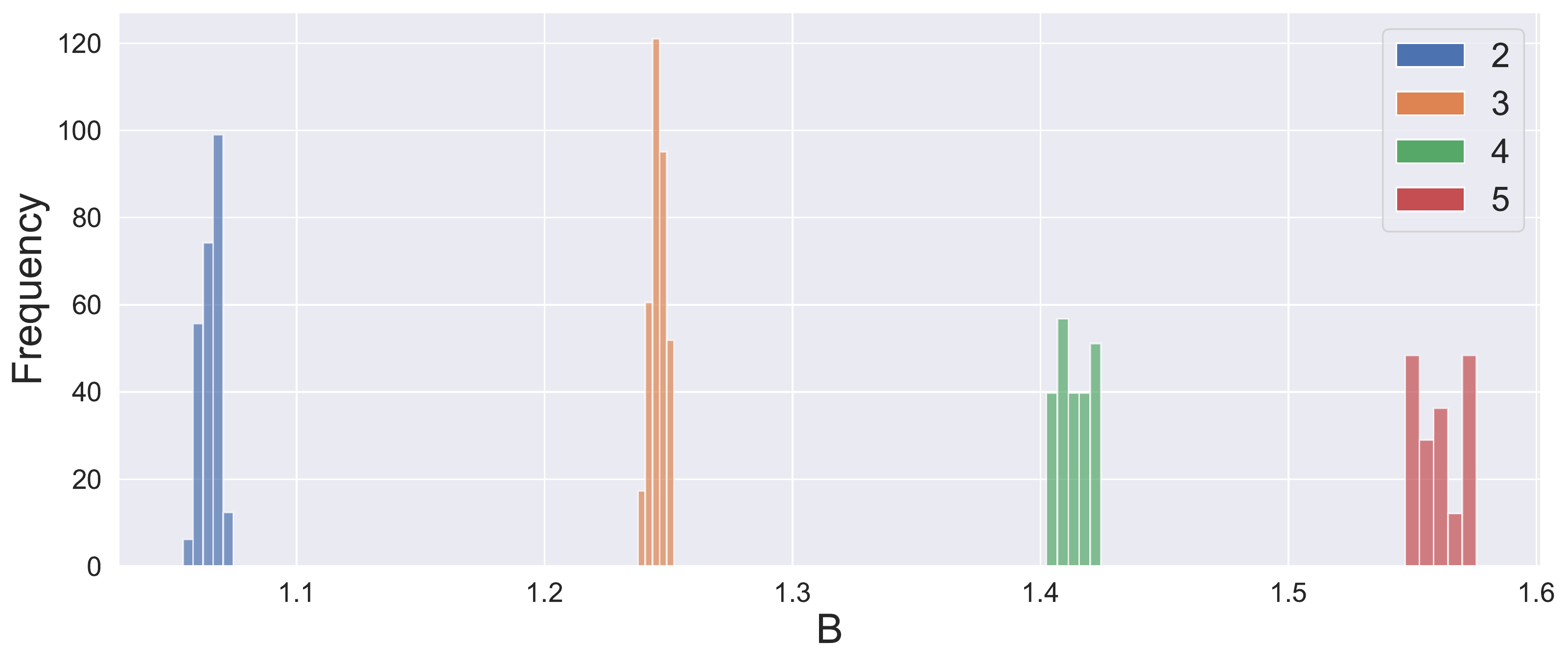}
     \end{subfigure}
     \begin{subfigure}{\textwidth}
         \includegraphics[width=0.95\textwidth]{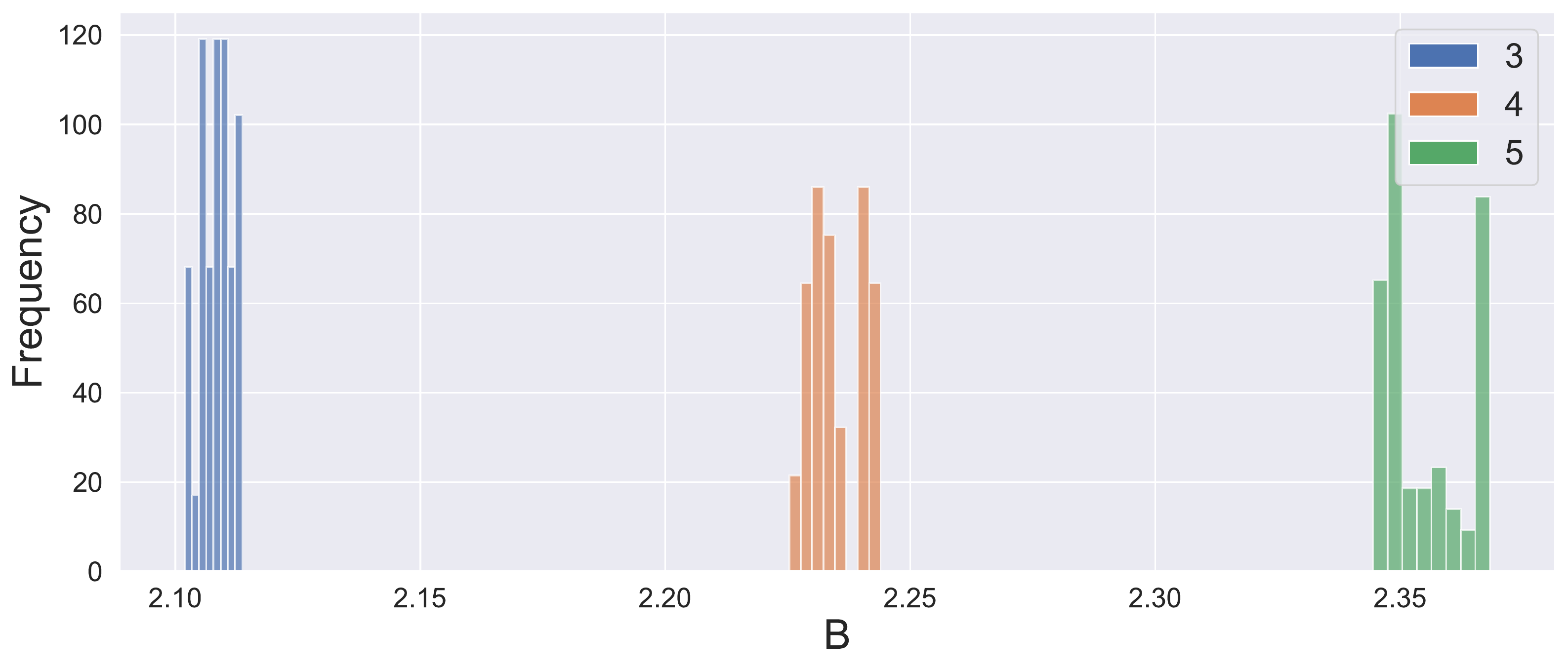}
     \end{subfigure}
     \begin{subfigure}{\textwidth}
           \includegraphics[width=0.95\textwidth]{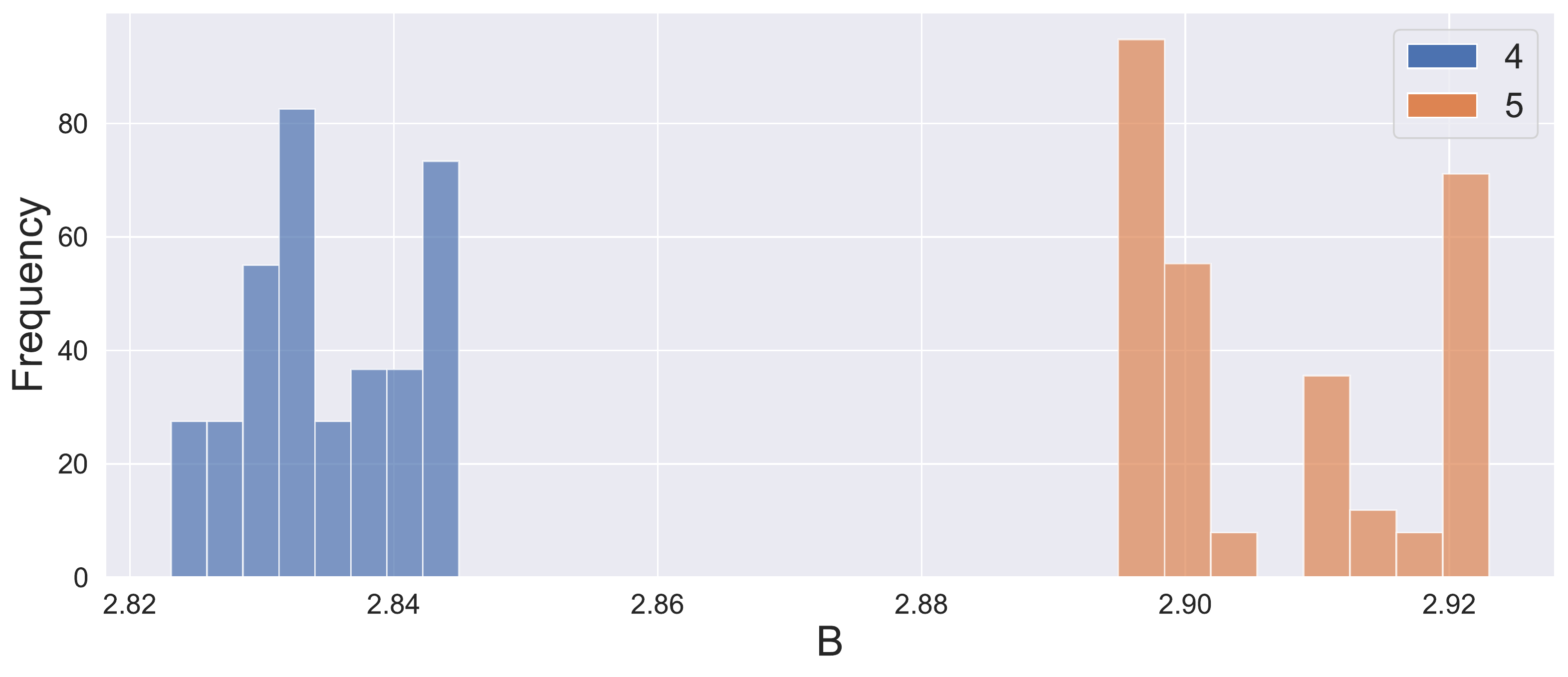}
     \end{subfigure}
       \caption{The distributions of the estimated $B$ values, for the \v Cech complex, across different models for $k=1$ (top), $k=2$ (middle), and $k=3$ (bottom).
    \label{fig:Alpha_H1_B}}
 \end{figure}

In Figures \ref{fig:Rips_H1_B} and \ref{fig:Alpha_H1_B}, we show the distributions of our estimates of $B$, across the various models, using sample size of $n=50,000$. As can be seen in the figures, the \v Cech complex estimates are much more concentrated than the Rips ones, and the distributions are more spread out in higher homological dimensions. We believe that this is due to slower convergence in these cases (in terms of the number of points).

\begin{figure}[h!]
     \centering
      \includegraphics[width=0.65\textwidth]{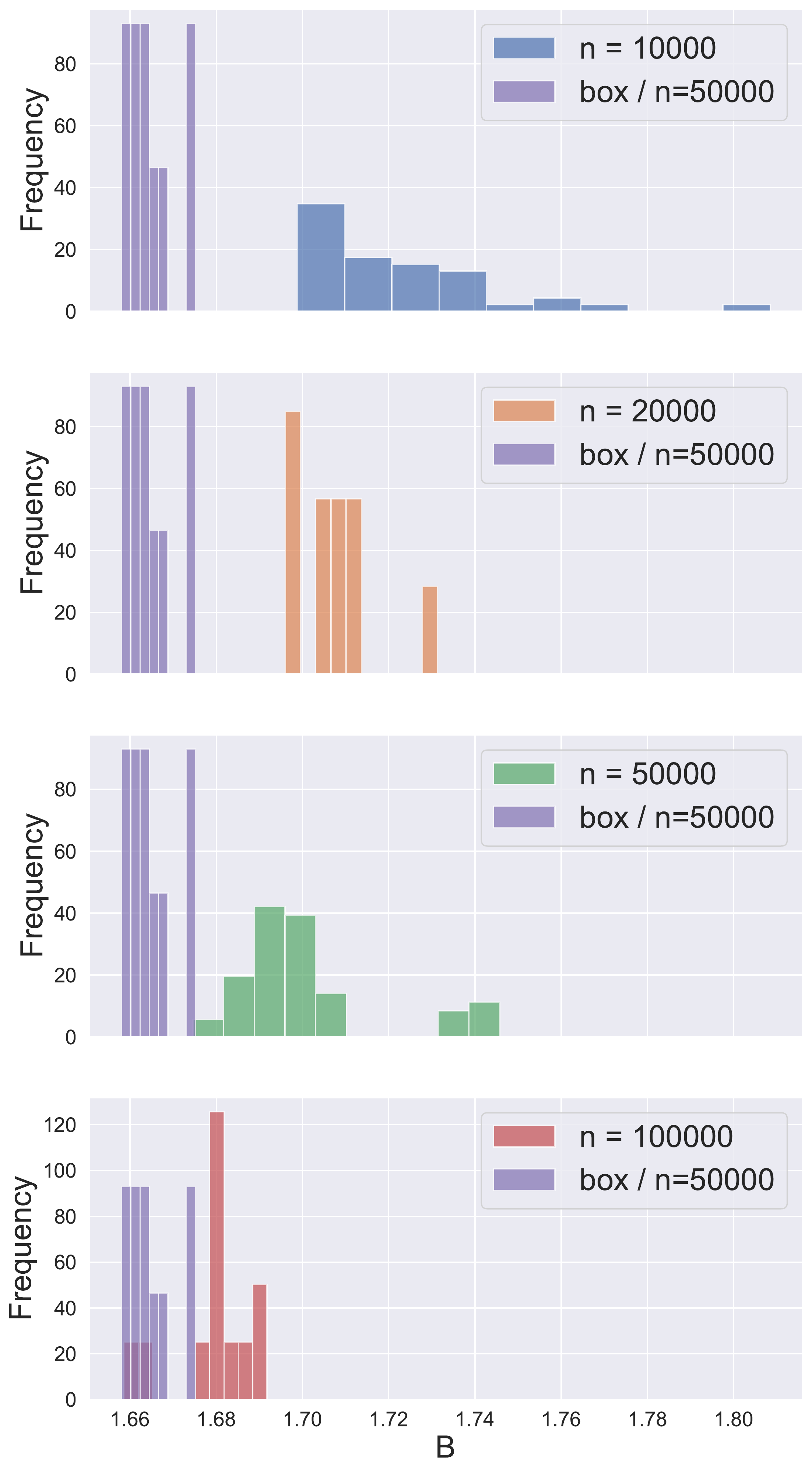}
       \caption{ The distributions of the estimated value of $B$, for different values of $n$.
       We compare the non-compact distributions in $\mathbb{R}^3$ (normal, Cauchy) for $n=10k, 20k, 50k, 100k$, with  the uniform distribution on $[0,1]^3$, $n=50k$ (in purple).
    \label{fig:N_dep}}
 \end{figure}

Finally, we wish to explore the connection between estimation variance and the sample size $n$.
In the case of Rips complex, the variance of the estimates of $B$ seems to be quite small in the compact-uniform case (e.g. on a box), compared to the variance emerging when sampling, for example, from the non-compact distribution (normal and Cauchy). In Figure~\ref{fig:N_dep}, we fix $d=3$ and $k=1$, and present the distribution for the box compared to the non-compact cases. The sample-size for the box was fixed $n=50,000$, and for the non-compact distribution we tested $n=10k 20k, 50k,100k$.
We observe that as $n$ increases the distribution of the non-compact cases moves closer to the box case. This serves as an evidence that in the limit all cases have indeed the same value of $B$ (supporting Conjecture \ref{con:iid}), it is just that some cases converge slower than others.

\end{document}